\documentclass[11pt,a4paper]{amsart}

\usepackage{amsmath,amssymb,amsthm,mathtools,mathrsfs}
\usepackage[T1]{fontenc}
\usepackage{lmodern}
\usepackage{microtype}
\usepackage[a4paper,margin=2.45cm,headheight=15pt]{geometry}
\usepackage{booktabs,array}
\usepackage{enumitem}
\usepackage{xcolor}
\usepackage{hyperref}
\usepackage[nameinlink,capitalise,noabbrev]{cleveref}

\definecolor{paperblue}{HTML}{173A5E}
\definecolor{papergold}{HTML}{8A5A16}
\hypersetup{
  colorlinks=true,
  linkcolor=black,
  citecolor=papergold,
  urlcolor=paperblue,
  pdftitle={Affine Vertex Operator Algebras and Correlation Functions on Elliptic Curves},
  pdfauthor={Linfang Hou, Si Li, Yongchang Zhu}
}

\allowdisplaybreaks[3]
\numberwithin{equation}{section}

\newtheorem{theorem}{Theorem}[section]
\newtheorem{maintheorem}{Theorem}

\newtheorem{proposition}[theorem]{Proposition}
\newtheorem{lemma}[theorem]{Lemma}
\newtheorem{corollary}[theorem]{Corollary}
\theoremstyle{definition}
\newtheorem{definition}[theorem]{Definition}

\theoremstyle{remark}
\newtheorem{remark}[theorem]{Remark}

\DeclareMathOperator*{\Res}{Res}
\DeclareMathOperator{\Tr}{Tr}
\DeclareMathOperator{\ad}{ad}

\newcommand{\CC}{\mathbb C}
\newcommand{\ZZ}{\mathbb Z}

\newcommand{\HH}{\mathbb H}
\newcommand{\g}{\mathfrak g}
\newcommand{\kk}{\Bbbk}
\newcommand{\one}{\mathbf 1}
\newcommand{\ii}{i}
\newcommand{\e}{e}
\newcommand{\dd}{\mathrm d}

\newcommand{\cP}{\mathcal P}

\newcommand{\scrD}{\mathscr D}
\newcommand{\zbf}{\boldsymbol z}
\newcommand{\ubf}{\boldsymbol u}
\newcommand{\mbf}{\boldsymbol m}

\title[Affine Current Traces on the Torus]{
Affine Vertex Operator Algebras and Correlation Functions on Elliptic Curves}

\author{Linfang Hou}
\address{L. Hou: Department of Mathematical Sciences, Tsinghua University, Beijing, China}
\email{houlf24@mails.tsinghua.edu.cn}

\author{Si Li}
\address{S. Li: Yau Mathematical Sciences Center, Tsinghua University, Beijing, China}
\email{sili@mail.tsinghua.edu.cn}

\author{Yongchang Zhu}
\address{Y. Zhu: Yau Mathematical Sciences Center and Beijing Institute of Mathematical Sciences and Applications,
Beijing, China}
\email{yc-zhu@mail.tsinghua.edu.cn}

\date{}

\begin{document}

\begin{abstract}
Let $V^{\kk}(\g)$ be the universal affine vertex algebra attached to a finite-dimensional complex simple Lie algebra at arbitrary level, and let $M=\bigoplus_{N=0}^{\infty}M_N$ be a $V^{\kk}(\g)$ -module.  We give an explicit formula for every genus-one trace of affine currents as a finite sum over cycle partitions.  We then show  the convergence of all fixed zero-mode trace series, hence absolute locally uniform convergence of the current traces and their meromorphic continuation on the elliptic configuration space.  The same conclusion is obtained for arbitrary affine-state insertions.  Finally, we compute the cases of at most four currents and the nodal limit $q\to 0$, in which the type-1 and type-2 elliptic cycle coefficients degenerate to rational genus-zero cycle coefficients and the module trace reduces to the ground-state space.
\end{abstract}
\maketitle

\tableofcontents

\section{Introduction}

Correlation functions play an important role in the theory of vertex operator algebras and conformal field theory.  In particular, genus-one $n$-point correlation functions encode important information about representations of vertex algebras. Their convergence, meromorphic continuation, and modular transformation properties are central topics in the theory.

In this paper, we study genus-one $n$-point correlation functions of currents associated with affine vertex algebras. Tsukada \cite{T} and Mason-Tuite \cite{MT} studied genus-one correlation functions for lattice vertex algebras, with the latter also obtaining explicit formulas for Heisenberg vertex algebras. In particular, the results of \cite{MT} implicitly provide formulas for the level-one affine vertex algebras of simply-laced type (ADE), through their lattice realizations. Bernard \cite{B} studied genus-one Ward identities for affine currents.  Dolan and Goddard \cite{DG07,DolanGoddard} investigated affine-current correlation functions and their applications to string theory.

In \cite{DolanGoddard}, Dolan and Goddard established a general formula for genus-one $n$-point correlation functions of affine currents at arbitrary level. Their formula is organized by permutations and their cycle decompositions, with each cycle contributing an elliptic function and an invariant tensor of the underlying Lie algebra. The elliptic functions are constructed from the Weierstrass sigma function, while the invariant tensors are determined by a systematic recursive procedure and are related to the characters of affine Lie algebra representations. Their work provides a general framework for computing torus current correlation functions and may be viewed as a genus-one analogue of the genus-zero current correlation formulas of Frenkel and the third author \cite{FrenkelZhu}.

The present paper develops a complementary approach to these correlation functions. One of our main results is a closed-form expression in terms of the Kronecker function and combinatorial sums over cycle decompositions, with coefficients given directly by finite-dimensional trace tensors. This formulation makes the dependence on the underlying finite-dimensional representation explicit and avoids the need to construct the invariant tensors recursively.

Let $\g$ be a finite-dimensional complex simple Lie algebra, and 
$h^\vee$ be its dual Coxeter number.  We fix the invariant form
\begin{equation}
\kappa(a,b)=\frac{1}{2h^\vee}\Tr_{\g}(\ad_a\ad_b),
\label{eq:normalized-form}
\end{equation}
so that every long root has squared length $2$.  The affine Kac-Moody Lie algebra is
\begin{equation*}
\widehat{\g}=\g\otimes\CC[t,t^{-1}]\oplus\CC K \oplus \CC d,
\end{equation*}
with
\begin{equation}
[j_m^a,j_n^b]=j_{m+n}^{[a,b]}+m\kappa(a,b)\delta_{m+n,0}K, \; \; [ d , j_m^a  ] = - m j_m^a , \; \; 
\qquad j_m^a:=a\otimes t^m.
\label{eq:affine-bracket-universal}
\end{equation}

Let $ M_{\rm top}$ be a finite dimensional $\g $-module, for a fixed level $ \kk  \in\CC$, 
we consider  $ M_{\rm top}$  as a module of $ 	\widehat{\g}_+ = \g\otimes\CC[t]  \oplus\CC K \oplus \CC d    $, where $ \g \otimes \CC [ t ] t $ 
acts as $0$, $K$ acts as $\kk$, and $ d $ acts as $ 0 $. We consider the induced representation 
\[  \widetilde{ M}_{\kk }  =   U ( \widehat{\g} )   \otimes_{ U ( \widehat{\g}_+ ) } M_{\rm top} . \]
We will consider any quotient representation $ M =  \widetilde{ M}_{ \kk  } /  I $ where $ I \cap M_{\rm top} = 0 $,
$M$  has a decomposition of $d$-eigenspaces:
\[   M   =  \bigoplus_{N=0}^{\infty}  M_N    \]
with $M_{0} = M_{\rm top}$ and $ M_N $ is spanned by elements
$   j_{-m_1}^{a_1} \cdots j_{-m_k}^{a_k}  v $, $v\in M_{\rm top}$ with $ m_i \geq 1 $ and $ \sum_{ i=1}^k  m _i = N $. 

In the case that $ M_{\rm top} = \CC$ is the trivial $ \g$-module,  we denote the induced representation by $V^{\kk} ( \g ) $, which has the structure of  a vertex algebra. This vertex algebra is known as the universal vertex algebra of level $\kk$ for $\g $. The generating vertex operators are
\begin{equation}
j^a(u)=Y(a,u)=\sum_{m\in\ZZ}j_m^a u^{-m-1},
\label{eq:plane-current}
\end{equation}
If the level $\kk\neq-h^\vee$, the vertex algebra $V^{\kk } ( \g ) $ has a Virasoro element $\omega $, it is given by  the Sugawara construction with 
central charge 
\begin{equation}
\omega=\frac{1}{2(\kk+h^\vee)}
\sum_\alpha j_{-1}^{x^\alpha}j_{-1}^{x_\alpha}\one,
\qquad
c=\frac{\kk\dim\g}{\kk+h^\vee},
\label{eq:sugawara}
\end{equation}
where $\{x_\alpha\}$ and $\{x^\alpha\}$ are dual bases. If $M_{\rm top}$ is irreducible,  the element $L_0 $ in the Virasoro algebra acts on $M_N$ as $N+h$, where $h$ is the scalar given by the $ \frac{1}{2(\kk+h^\vee)}$ times the Casimir element on $ M_{\rm top}$.

Throughout this paper, we omit the vacuum-energy factor from the trace
normalization: we will always use
$\Tr_M(Aq^d)$ rather than the customary physical trace
$\Tr_M(Aq^{L_0-c/24})$ even when the Virasoro element is available. The latter is obtained by multiplying the former
by $q^{h-c/24}$ $(L_0|_{M_0}=h\,\mathrm{id}_{M_0})$.

For every finite dimensional $\g $-module $M_{\rm top}$, every quotient  $M =   \widetilde{ M}_{ \kk  } /  I $ carries the structure of a module over the vertex algebra $ V^{\kk} ( \g )$. The purpose of this work is to study the $n$-point trace 
\begin{equation} \label{n-trace}
( 2\pi i )^n   u_1 \cdots u_n   \Tr_M  j^{a_1} ( u_1 ) \cdots j^{a_n} ( u_n ) q^d .   
\end{equation}
Our first main result is

\begin{maintheorem}\label{theoremA}
	The $n$-point trace (\ref{n-trace}) converges in the region $  0 < |q| < |u_n| < \dots < |u_1| < 1 $ and the limit admits a meromorphic continuation to the region $ 0 <|q|<1 ,  u_1 , \dots , u_n \in \CC^* $. Moreover the meromorphic continuation is invariant under the transformation $ u_i \mapsto q^k u_i $ for $ k \in \ZZ $ and $1\leq i \leq n$.
\end{maintheorem}

If a vertex operator algebra satisfies the $C_2$-condition, the convergence of an $n$-point trace is proved in \cite{Zhu}.
The proof is roughly divided into two steps. Step 1 reduces the problem to proving the convergence of the 1-point trace. Step 2 shows that the 
$C_2$-condition implies that the 1-point trace satisfies a differential equation, which in turn yields its convergence. Subsequently, it was noticed in \cite{AK} that the same method can be used to establish convergence if the 
$C_2$-condition is replaced by the weaker condition
\[ {\rm dim}  V / \mathrm{Span} ( a_{(-2)} b, a_{(0)} b \, | \, a , b \in V )  < \infty .\]        
In our case, the vertex algebra is $V^{\kk}(\frak{g})$, and the above finiteness conditions are not satisfied in general—not even if $V^{\kk}(\frak{g})$ is replaced by its irreducible quotient $L_{\kk}(\frak{g})$. We refer to \cite{Lin} for a comprehensive discussion of the 
$C_2$-condition for the vertex algebra $L_{\kk}(\frak{g})$.

Writing
\begin{equation*}
q=\e^{2\pi\ii\tau},\qquad
u=\e^{2\pi\ii z},\qquad \tau\in\HH.
\end{equation*}
We define
\begin{equation}\label{n-function}
   S_{n,M}^{a_1 , \dots , a_n }  (  u_1 , \dots , u_n ; q    ) =  ( 2\pi i )^n   u_1 \cdots u_n   \Tr_M  j^{a_1} ( u_1 ) \cdots j^{a_n} ( u_n ) q^d     
\end{equation}   
initially in the region $  0 < |q| < |u_n| < \dots < |u_1| < 1 $ and extend $S_M$ to the region $ |q|<1 ,  u_1 , \dots , u_n \in \CC^* $. 
The functions $ S_M$ are elliptic functions with periods $ \ZZ +\ZZ \tau $ and possible poles only where $z_i-z_j\in\ZZ+\tau\ZZ$ for $i\neq j$. 
When no confusion occurs, we omit the subscript $M$ and write the function
$S_M $ as $S$. 

Our second main result expresses the $n$-point trace function $S_M$  as a finite sum over partitions of $ [ n ]= \{1,\ldots,n\}$ into cyclically ordered subsets. Each cycle is assigned one of two types of weights, defined in terms of the Kronecker function and certain trace functions associated with the zero modes of the affine Lie algebra. We now introduce the notation needed to state this formula.
Recall the odd Jacobi theta function is
\begin{equation}
\theta(z,\tau)=-\ii\sum_{m\in\ZZ}(-1)^m
q^{\frac12(m+\frac12)^2}\e^{2\pi\ii(m+\frac12)z},
\label{eq:odd-theta-intro}
\end{equation}
We introduce the Kronecker function
\begin{equation}
K(z,\lambda;\tau)
:=\frac{\theta'(0,\tau)\theta(z+\lambda,\tau)}
{\theta(z,\tau)\theta(\lambda,\tau)}.
\label{eq:Kronecker-intro}
\end{equation}
It obeys
\begin{equation}
K(z+1,\lambda)=K(z,\lambda),\qquad
K(z+\tau,\lambda)=\e^{-2\pi\ii\lambda}K(z,\lambda),
\label{eq:Kronecker-periods-intro}
\end{equation}
and $K(z,\lambda)=z^{-1}+O(1)$ as $z\to0$.  See \cite{Fay} for more details.
For a cyclically ordered set $ C = ( i_1 i_2 \dots i_r)$, where $\{ i_1 , \dots , i_r \}$ is a subset of $ [n]= \{ 1 , \dots , n \}$, we introduce the elliptic function in $z_i$ associated with $C$ 
 by  
\begin{equation*}
H_r ( C; z_{i_1} , \dots , z_{i_r} ; \lambda  ) =\lambda^r\prod_{s=1}^r
K(z_{i_s}-z_{i_{s+1}},\lambda;\tau),
\qquad i_{r+1}=i_1,
\end{equation*}
We also denote this by  $H_C ( \lambda ) $. 
We introduce the \emph{type-1 elliptic cycle coefficient} $\Psi_C$ and the \emph{type-2 elliptic cycle coefficient} $\Phi_C$, functions of
$z_{i_1},\ldots,z_{i_r},\tau$, by
\begin{align}
\Psi_r (C , z_{i_1} , \dots , z_{i_r} , \tau )  &= \Psi_r(C ) 
=2\pi\ii[\lambda^{r-1}]
\bigl(\e^{2\pi\ii\lambda}H_r(C;\lambda)\bigr), \label{1.10}\\
\Phi_r(C, z_{i_1} , \dots , z_{i_r} , \tau )&=\Phi_r(C) = -[\lambda^r]\left(
\lambda\frac{\theta'(\lambda,\tau)}{\theta(\lambda,\tau)}
H_r(C;\lambda)\right). \label{1.11}
\end{align}
where for a meromorphic function $f$ of variables $\lambda , z_1 , \dots , z_r, \tau $ 
$ [ \lambda^k ] f ( \lambda , z_{1} , \dots , z_{r})$
denotes the coefficient of $ \lambda^k $ in the Laurent expansion of $f$ in the variable $\lambda $. So this coefficient is a meromorphic function of $z_1 , \dots , z_r $.  See Section 3 for more properties of these functions.

For each cyclic set $C$ as above, we introduce the following two types of weights.
They are 
\begin{equation*}
W_1(C)=t_M(j_0^{a_{i_1}}\cdots j_0^{a_{i_r}};q)\Psi_C ( z_{i_1 } , \dots , z_{i_r } ; \tau )
\qquad
W_2(C)=\kk\kappa(a_{i_1},\ldots,a_{i_r})\Phi_C ( z_{i_1 } , \dots , z_{i_r } ; \tau )
\end{equation*}
where
\begin{equation*}
t_M(j_0^{a_{i_1}}\cdots j_0^{a_{i_r}};q)  :=\sum_{N =0}^\infty \Tr_{M_{N}}(j_0^{a_{i_1}}\cdots j_0^{a_{i_r}} )q^N,
\qquad
\kappa(a_1,\ldots,a_r):=\frac{1}{2h^\vee}
\Tr_\g(\ad_{a_1}\cdots\ad_{a_r}).
\end{equation*}
Let $\cP_n^\star$ be  the set of  cycle partitions of the set $[n]= \{ 1 , \dots , n \}$, where each cycle has length at least two and is labeled type 1 or type 2, with at most one cycle of type 1.  For each such decomposition $\alpha $ of $[n]$:
\[  \alpha = \{  C_1 ,  \dots , C_k \}   \]
we define 
\[  W ( \alpha ) = \begin{cases}
\displaystyle W_1(C_\star)
\prod_{\substack{C\in\alpha\\C\neq C_\star}}W_2(C),
&C_\star \text{ is the unique type-1 cycle},\\[2ex]
\displaystyle t_M(1;q)\prod_{C\in\alpha}W_2(C),
& \text{all the cycles are of type 2 }.
\end{cases}                            \]

Now we can state our main theorem

\begin{maintheorem}	\label{theoremB}
	For every $n\geq0$ and $a_1,\ldots,a_n\in\g$, the correlation function of currents is given by
	\begin{equation*}
	S_{n, M}^{ a_1 , \dots , a_n} ( e^{2\pi i z_1}  , \dots ,e^{2\pi i z_n}  ; e^{2\pi i \tau } ) 
	=\sum_{\alpha\in\cP_n^\star}  W ( \alpha )
	\end{equation*}
\end{maintheorem}

This formula is similar to a formula about genus zero correlation function for vertex algebra $ V^\kk ( \g )$. The more precise relation is that $ \lim_{ \tau \to i \infty }	S_{M} ( a_1 , \dots , a_n; z_1 , \dots , z_n ; \tau ) $ reduces to the formula in \cite{FrenkelZhu}.

\

The Proof of Theorem A  and Theorem B can be described as follows.  We denote the right-hand side of the identity in Theorem B as $ G_{n ,M}^{a_1 , \dots , a_n}$, we prove two sides are equal as a formal series by  comparing the coefficients of power series on each side. The coefficients of the left hand sides are denoted by $ T^S_M ( w , q )$, where $w $ is a word in $ \widehat{\frak g}$ (which is called finite mode words),
 the coefficient of the right hand is called $T^G_M ( w , q ) $. Our strategy is to prove 
 \[     T^S_M ( w , q ) = T^G_M ( w , q ) \]
To this end, we show that the system of coefficients of $ T^S_M ( w , q )$ satisfies the  relations of stated in Proposition 4.1. We then prove that
 the system of coefficients of  $T^G_M ( w , q )$ satisfies the same relations, this is the most technical part of this works, and we complete the proof in 
Proposition 4.3.  Since these relations imply that 
$T_M^S ( w , q ) $ and $ T_M^G ( w , q )$
are determined by the corresponding coefficients of length $1$
, and since it is clear that $T_M^S ( w , q ) = T_M^G ( w , q )$ for $w$ of length $1$,
 we have proved the desired identity (Theorem 4.4).
  The equality of the $ n$-variable functions 
 $S_{n, M }$ and $G_{n, M }$, together with the convergence results for the 
 $1$-point trace function proved in Lemma 5.5, implies Theorem A, and Theorem B then follows automatically. 
 
 Using Theorems A and B, we establish our third main result,
  Theorem \ref{thm:zhu-recursion-general-convergence}, it proves  the convergence and meromorphc continuation of $n$-point trace of arbitrary vertex operators 
 \begin{equation}
 S_{n,M}^{v_1,\ldots,v_n}(\ubf;q)
 :=\left(\prod_{i=1}^n(2\pi\ii u_i)^{|v_i|}\right)
 \Tr_M \left(Y_M(v_1,u_1)\cdots Y_M(v_n,u_n)q^d\right).
 \end{equation}

 Finally, we discuss some potential applications of our main results. For applications of affine-current correlation functions in mathematical physics, we refer to \cite{DolanGoddard} and the references therein. Here we focus on applications to the representation theory of affine vertex algebras.
 
 It is known that the characters of ordinary modules over a quasi-lisse vertex operator algebra satisfy modular linear differential equations under suitable hypotheses \cite{AK}. Quasi-lisse simple affine vertex algebras $L_k(\frak g)$ have recently been classified in \cite{ADS}. By definition, the associated variety of a quasi-lisse vertex algebra has only finitely many symplectic leaves. This raises the question of whether certain symplectic leaves can be related to modules over $L_k(\frak g)$ whose lowest conformal-weight spaces $M_{\rm top}$ are infinite-dimensional.
 
 In this setting, an interesting problem is to extend our formulas to such modules. In particular, Theorem B suggests that the construction of genus-one correlation functions requires a suitable interpretation of the zero-mode traces
 \[
 \operatorname{Tr}_M \left(
 j_0^{a_1}\cdots j_0^{a_n}q^d
 \right).
 \]
 When $M_{\rm top}$ is infinite-dimensional, these traces need not converge or even be defined in the usual sense. It would therefore be interesting to investigate possible regularizations or alternative trace functionals that would allow our formulas to be extended to this setting.

\

The paper is organized as follows. In Section 2, we recall some basic facts about elliptic functions and introduce the multi-variable elliptic functions 
$\Psi (C) $ and $\Phi ( C)$ 
that appear in our main formula. In Section 3, we construct a family of elliptic functions, which we call the graph sum of $n$-point correlators  associated with a representation 
$M$ of the affine Lie algebra, which will be identified with the 
$n$-point trace functions. In Section 4, we establish this identification at the level of formal power series. In Section 5, we use these formal identities to prove the convergence and meromorphic continuation of the 
$n$-point trace functions, thereby completing the proofs of our main theorems. In Section 6, we discuss our formula for 
$n$-point functions in the cases $n \leq 4 $.
 In Section 7, we consider the limit $\tau \to i \infty $
and show in Theorem 7.3 that the limit of our formula agrees with the genus-zero 
$n$-point correlation functions given in \cite{FrenkelZhu}. Appendix A contains some facts about Bell and Bernoulli polynomials used in the proofs of our main theorems.

\

\subsection*{Acknowledgements}
S. L. is supported by NSFC No.12571068. Part of the exposition in Section 4 was prepared with AI assistance, under the author's supervision and verification.

\section{Elliptic Functions and Cycle Coefficients}
\label{sec:elliptic-kernels}

In this section, we recall some basic facts about elliptic functions, including the quasiperiodic Weierstrass zeta function $\wp_1$, the Weierstrass $\wp$-function $\wp_2$, the odd Jacobi theta function $\theta$, and the Kronecker function $K$. We then introduce the type-1 elliptic cycle coefficient $\Psi(C)$ and the type-2 elliptic cycle coefficient $\Phi(C)$ associated to a cycle $C$, which will appear in our explicit formula for the $n$-point trace functions.

\subsection{Eisenstein, theta, and Weierstrass functions}

Let $\Lambda_\tau=\ZZ+\tau\ZZ$ and
$E_\tau=\CC/\Lambda_\tau$.  For $r\geq2$, define
\begin{equation*}
 G_{2r}(\tau)=
 \sum_{(m,n)\in\ZZ^2\setminus\{(0,0)\}}
 \frac{1}{(m+n\tau)^{2r}},
\end{equation*}
and use Eisenstein summation for
\begin{equation*}
 G_2(\tau)=\sum_{n\neq0}\frac1{n^2}
 +\sum_{m\neq0}\sum_{n\in\ZZ}\frac1{(m\tau+n)^2}.
\end{equation*}
The odd theta function is
\begin{equation}
 \theta(z,\tau)=-\ii\sum_{m\in\ZZ}(-1)^m
 q^{\frac12(m+\frac12)^2}\e^{2\pi\ii(m+\frac12)z},
 \label{eq:odd-theta}  
\end{equation}
and satisfies
\begin{equation*}
 \theta(z+1,\tau)=-\theta(z,\tau),\qquad
 \theta(z+\tau,\tau)=-q^{-1/2}\e^{-2\pi\ii z}\theta(z,\tau).
\end{equation*}
These conventions of theta function here differ from that with characteristic $(\frac12,\frac12)$ in \cite[Chapter~I]{Mumford} by a minus sign.
We have the product formula
\begin{equation}
\theta(z,\tau)
= i q^{1/8} e^{-\pi i z}
(1-e^{2\pi i z})
\prod_{n=1}^{\infty}
(1-q^n)
(1-q^n e^{2\pi i z})
(1-q^n e^{-2\pi i z}).
\end{equation}

Define the Weierstrass sigma function by
\begin{equation}
 \sigma(z,\tau)=\frac{\theta(z,\tau)}{\theta'(0,\tau)}
 \exp\!\left(\frac{G_2(\tau)z^2}{2}\right),
 \label{eq:sigma-theta}
\end{equation}
and set
\begin{equation}
 \wp_1(z,\tau)=\partial_z\log\sigma(z,\tau),
 \qquad
 \wp_2(z,\tau)=-\partial_z\wp_1(z,\tau).
 \label{eq:wp1-wp2}
\end{equation}
Thus $\wp_1$ is the Weierstrass zeta function and $\wp_2$ is the usual
Weierstrass $\wp$ function.  From \eqref{eq:sigma-theta},
\begin{align}
 \wp_1(z,\tau)
 &=\frac{\theta'(z,\tau)}{\theta(z,\tau)}+G_2(\tau)z,
 \label{eq:wp1-theta}\\
 \wp_2(z,\tau)+G_2(\tau)
 &=-\partial_z^2\log\theta(z,\tau).
 \label{eq:wp2-theta}
\end{align}
The local expansions are
\begin{align}
 \wp_1(z)&=\frac1z-\sum_{r\geq2}G_{2r}z^{2r-1},
 \label{eq:wp1-expansion}\\
 \wp_2(z)&=\frac1{z^2}
 +\sum_{r\geq1}(2r+1)G_{2r+2}z^{2r},
 \label{eq:wp2-expansion}\\
 \log\frac{\sigma(z)}z
 &=-\sum_{r\geq2}\frac{G_{2r}}{2r}z^{2r}.
 \label{eq:sigma-expansion}
\end{align}
These standard formulas may be found in
\cite{WhittakerWatson} and \cite{Fay}.

The functions appearing in the recurrent formula \eqref{prop:zhu-recursion-recalled} are
\begin{equation*}
 \wp_{m+1}(z,\tau)
 :=\frac{(-1)^m}{m!}\,\partial_z^m\wp_1(z,\tau).\qquad m\geq 0.
\end{equation*}

We will frequently meet 
\begin{equation}
 P_1(z,\tau):=-\wp_1(z,\tau)+G_2(\tau)z-\pi\ii,
 \qquad
 P_2(z,\tau):=\partial_zP_1(z,\tau).
 \label{eq:P1P2}
\end{equation}
Equivalently,
\begin{equation*}
 P_1(z)=-\frac{\theta'(z)}{\theta(z)}-\pi\ii,
 \qquad P_2(z)=\wp_2(z)+G_2.
\end{equation*}
Their local singular parts are
\begin{equation}
 -P_1(z)=\frac1z+\pi\ii+O(z),
 \qquad P_2(z)=\frac1{z^2}+G_2+O(z^2),
 \label{eq:P-principal-parts}
\end{equation}
and, with $u=\e^{2\pi\ii z}$,
\begin{align}
 P_1(u;q)&=(2\pi\ii)\sum_{m\geq1}
 \left(\frac{u^m}{1-q^m}-\frac{u^{-m}q^m}{1-q^m}\right),
 \label{eq:P1-Fourier}\\
 P_2(u;q)&=(2\pi\ii)^2\sum_{m\geq1}
 \left(\frac{m u^m}{1-q^m}+\frac{m u^{-m}q^m}{1-q^m}\right).
 \label{eq:P2-Fourier}
\end{align}
The series converge absolutely and locally uniformly for $|q|<|u|<1$; they are the $P$-functions introduced in \cite{Zhu}.

\subsection{The Kronecker function}

Define
\begin{equation}
 K(z,\lambda;\tau)
 :=\frac{\theta'(0,\tau)\theta(z+\lambda,\tau)}
 {\theta(z,\tau)\theta(\lambda,\tau)}.
 \label{eq:Kronecker}
\end{equation}
It obeys
\begin{equation}
 K(z+1,\lambda)=K(z,\lambda),\qquad
 K(z+\tau,\lambda)=\e^{-2\pi\ii\lambda}K(z,\lambda),
 \label{eq:Kronecker-periods}
\end{equation}
and $K(z,\lambda)=z^{-1}+O(1)$ as $z\to0$.  See \cite{Fay} for more details. The cyclic product below is the generating function used in \cite{DolanGoddard}.

A cycle is a finite subset of $[n]=\{1,\ldots,n\}$ of cardinality at least two, equipped with a cyclic order.  For
$C=(i_1\,i_2\cdots i_r)$, put
\begin{equation*}
 \mu_s=z_{i_s}-z_{i_{s+1}},\qquad i_{r+1}=i_1,
 \qquad \sum_{s=1}^r\mu_s=0,
\end{equation*}
and define
\begin{equation}
 H_r(C;\lambda):=\lambda^r\prod_{s=1}^rK(\mu_s,\lambda;\tau)
 =\sum_{m\geq0}H_{r,m}(C;\tau)\lambda^m.
 \label{eq:H-definition}
\end{equation}
Using \eqref{eq:sigma-theta} on every edge gives
\begin{equation}
 H_r(C;\lambda)=\frac{\lambda^r}{\sigma(\lambda)^r}
 \prod_{s=1}^r\frac{\sigma(\mu_s+\lambda)}{\sigma(\mu_s)}.
 \label{eq:H-sigma}
\end{equation}

If the edge $\mu_r$ shrinks, then for $r>2$,
\begin{equation}
 \Res_{\mu_r=0}H_r(C;\lambda)\,\dd\mu_r
 =\lambda H_{r-1}(C/(i_ri_1);\lambda),
 \label{eq:H-collision}
\end{equation}
where the two adjacent vertices are fused.  Hence
\begin{equation}
 \Res_{\mu_r=0}H_{r,m}(C)\,\dd\mu_r=H_{r-1,m-1}(C/(i_ri_1)).
 \label{eq:H-coefficient-collision}
\end{equation}
Moreover, $\lambda^{-r}H_r=\prod_sK(\mu_s,\lambda)$ is elliptic in
$\lambda$.  Its only pole in a fundamental parallelogram is at
$\lambda=0$, and its residue is $H_{r,r-1}$.  The residue theorem gives
\begin{equation}
 H_{r,r-1}=0\qquad(r\geq2).
 \label{eq:H-vanishing}
\end{equation}

For explicit coefficients, write
\begin{equation*}
 \log H_r(C;\lambda)=\sum_{s\geq1}L_s(C)\frac{\lambda^s}{s!}.
\end{equation*}
Equations \eqref{eq:wp1-expansion}--\eqref{eq:sigma-expansion} give
\begin{equation}
 L_s(C)=\sum_{j=1}^r\wp_1^{(s-1)}(\mu_j)
 +\begin{cases}
 r(s-1)!G_s,&s\geq4\text{ even},\\
 0,&\text{otherwise}.
 \end{cases}
 \label{eq:L-coefficients}
\end{equation}
If $B_m$ is the complete exponential Bell polynomial of
\cref{app:bell}, then
\begin{equation}
 H_{r,m}=\frac1{m!}B_m(L_1,\ldots,L_m)
 \label{eq:H-Bell}
\end{equation}
gives a formula for exact computation.

\subsection{Type-1 and type-2 elliptic cycle coefficients}

For a directed cycle $C$ of length $r\geq2$, as in (\ref{1.10}) (\ref{1.11}), we define the
\emph{type-1 elliptic cycle coefficient} $\Psi_r(C)$ and the
\emph{type-2 elliptic cycle coefficient} $\Phi_r(C)$ by
\begin{align}
 \Psi_r(C)
 &:=2\pi\ii[\lambda^{r-1}]
 \left(\e^{2\pi\ii\lambda}H_r(C;\lambda)\right),
 \label{eq:Psi-definition}\\
 \Phi_r(C)
 &:=-[\lambda^r]\left(
 \lambda\frac{\theta'(\lambda,\tau)}{\theta(\lambda,\tau)}
 H_r(C;\lambda)\right).
 \label{eq:Phi-definition}
\end{align}
Here square brackets denote the extraction of the coefficient of the corresponding term in the expansion. We also set
\begin{equation*}
 \Psi_1:=2\pi\ii,\qquad \Phi_1:=0.
\end{equation*}
In terms of \eqref{eq:H-definition},
\begin{align}
 \Psi_r(C)&=\sum_{m=0}^{r-1}
 \frac{(2\pi\ii)^{r-m}}{(r-m-1)!}H_{r,m}(C),
 \label{eq:Psi-H}\\
 \Phi_r(C)&=-H_{r,r}(C)
 +\sum_{s=1}^{\lfloor r/2\rfloor}G_{2s}H_{r,r-2s}(C).
 \label{eq:Phi-H}
\end{align}
The final term in \eqref{eq:Psi-H} vanishes by
\eqref{eq:H-vanishing}.

\section{The Graph Sum $n$-Point Correlators}
\label{sec:graph-sum}

In this section, we introduce the graph sum $n$-point correlators of affine currents $j^{a_1} , \dots , j^{a_n}$ and a $\hat {\frak g}$-module as in Section 1.

For $r\geq2$, extend the Killing form to a series of maps $\mathfrak{g}^{\otimes r}\to \CC$
\begin{equation}
 \kappa(a_1,\ldots,a_r):=\frac1{2h^\vee}
 \Tr_{\g}(\ad_{a_1}\cdots\ad_{a_r}).
 \label{eq:extended-form}
\end{equation}
Clearly for $r\geq 3$ we have 
\begin{align}
 \kappa(a_1,a_2,\ldots,a_r)&=\kappa(a_2,\ldots,a_r,a_1),
 \label{eq:extended-cyclic}\\
 \kappa(a,b,a_3,\ldots,a_r)-\kappa(b,a,a_3,\ldots,a_r)
 &=\kappa([a,b],a_3,\ldots,a_r).
 \label{eq:extended-commutator}
\end{align}

 Let $C=(i_1\,i_2\cdots i_r)$ be a cycle.  We will use two types of cycles.

\begin{enumerate}
\item
The type--1 cycle records the trace information of the zero mode. If we set the formal power series of zero mode trace as
\begin{equation}
 t_M(j_0^{a_1}\cdots j_0^{a_r};q):=
 \sum_{N\geq0}\Tr_{M_{ N}}
 (j_0^{a_1}\cdots j_0^{a_r})q^N,
 \qquad t_M(1;q):=\sum_{N\geq0}\dim M_{N}q^N,
 \label{eq:zero-mode-series-early}
\end{equation}

then the corresponding weight associated to a type-1 cycle is
\[
W_1(C)=
t_M(j_0^{a_{i_1}}\cdots j_0^{a_{i_r}};q) \Psi_r(C)
\]
\item
The type--2 cycle records the information of the level $\kk$ and extended Killing form, and the corresponding weight is
\[
W_2(C)=
\kk\kappa(a_{i_1},\ldots,a_{i_r})\Phi_r(C)
\]
\end{enumerate}

\begin{figure}[ht]
\centering
\begin{minipage}[t]{0.47\textwidth}
\centering
\setlength{\unitlength}{0.75pt}
\begin{picture}(205,105)
  \linethickness{0.8pt}
  \put(30,75){\circle{22}}
  \put(30,75){\makebox(0,0){$i_1$}}
  \put(150,75){\circle{22}}
  \put(150,75){\makebox(0,0){$i_2$}}
  \put(150,18){\circle{22}}
  \put(150,18){\makebox(0,0){$i_3$}}
  \put(30,18){\circle{22}}
  \put(30,18){\makebox(0,0){$i_r$}}
  \put(42,75){\vector(1,0){96}}
  \put(150,63){\vector(0,-1){32}}
  \put(138,18){\vector(-1,0){32}}
  \put(90,18){\makebox(0,0){$\cdots$}}
  \put(74,18){\vector(-1,0){32}}
  \put(30,30){\vector(0,1){32}}
  \put(90,46){\makebox(0,0){\color{black}\bfseries type 1}}
  \put(90,2){\makebox(0,0){\small one distinguished cycle}}
\end{picture}
\[
W_1(C)=
t_M(j_0^{a_{i_1}}\cdots j_0^{a_{i_r}};q)\Psi_r(C)
\]
\end{minipage}
\hfill
\begin{minipage}[t]{0.47\textwidth}
\centering
\setlength{\unitlength}{0.75pt}
\begin{picture}(205,105)
  \linethickness{0.8pt}
  \put(30,75){\circle{22}}
  \put(30,75){\makebox(0,0){$i_1$}}
  \put(150,75){\circle{22}}
  \put(150,75){\makebox(0,0){$i_2$}}
  \put(150,18){\circle{22}}
  \put(150,18){\makebox(0,0){$i_3$}}
  \put(30,18){\circle{22}}
  \put(30,18){\makebox(0,0){$i_r$}}
  \put(42,75){\vector(1,0){96}}
  \put(150,63){\vector(0,-1){32}}
  \put(138,18){\vector(-1,0){32}}
  \put(90,18){\makebox(0,0){$\cdots$}}
  \put(74,18){\vector(-1,0){32}}
  \put(30,30){\vector(0,1){32}}
  \put(90,46){\makebox(0,0){\color{black}\bfseries type 2}}
  \put(90,2){\makebox(0,0){\small a closed central cycle}}
\end{picture}
\[
W_2(C)=
\kk\kappa(a_{i_1},\ldots,a_{i_r})\Phi_r(C)
\]
\end{minipage}
\caption{The two type cycle weights.}
\label{fig:cycle-types}
\end{figure}
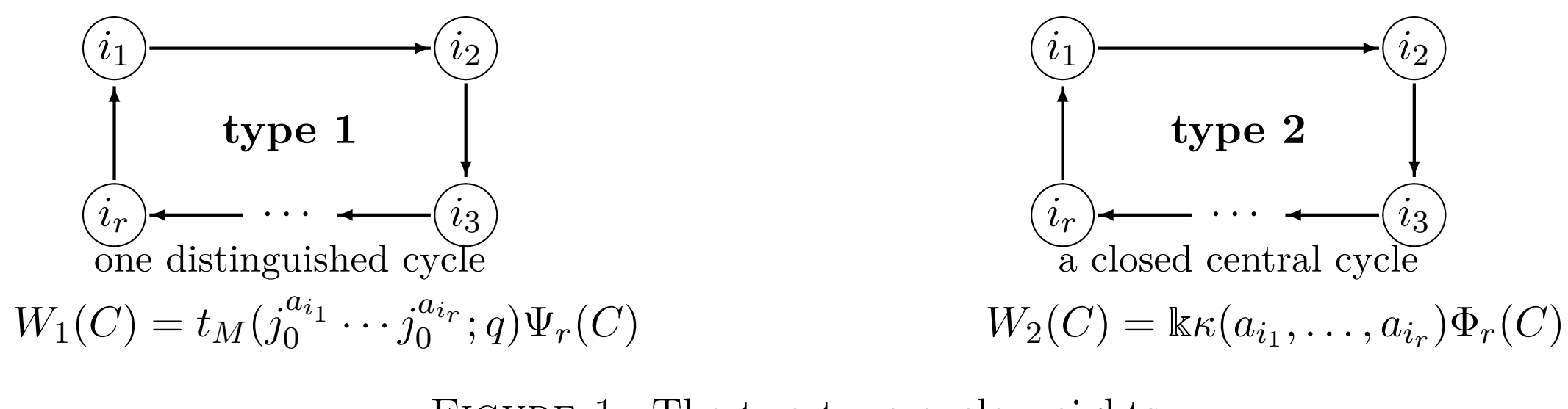

Note that both $t_M(j_0^{a_{i_1}}\cdots j_0^{a_{i_r}};q)$ and $\kk\kappa(a_{i_1},\ldots,a_{i_r})$ are cyclically symmetric by definition.

Let $[n]=\{1,\ldots,n\}$.  A \textbf{cycle partition} $\alpha$ is a partition of $[n]$ into disjoint cycles and each is labeled type--1 or type--2. 

Let $\cP_n^\star$ denote the set of all cyclic partitions with \textbf{at most $1$} cycle of type--1.

If $\alpha\in \cP_n^\star$ has a
type--1 cycle $C_\star$, define the weight of the partition $\alpha$
\begin{equation*}
  W(\alpha)
  :=W_1(C_\star)
    \prod_{\substack{C\in\alpha\\C\neq C_\star}}W_2(C).
\end{equation*}
If no type--1 cycle is present, define
\begin{equation*}
  W(\alpha)
  :=t_M(1;q)\prod_{C\in\alpha}W_2(C).
\end{equation*}

So a general term is of the form
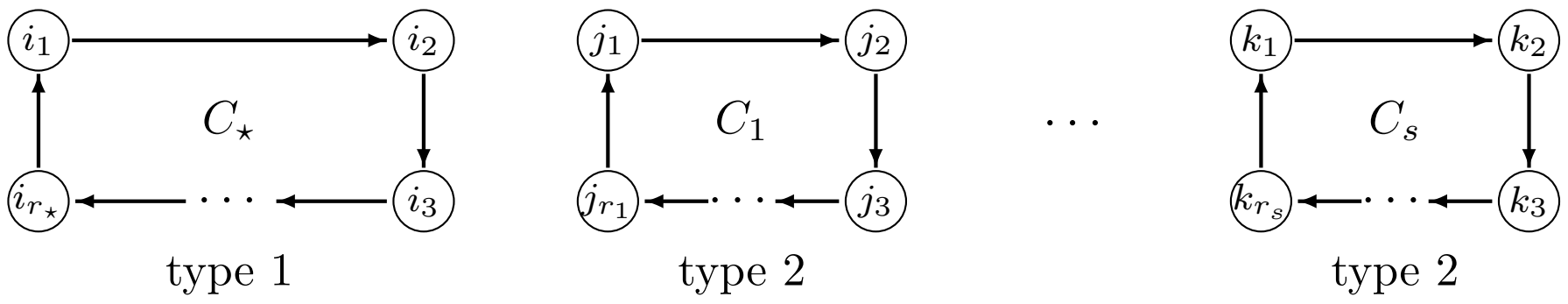
\begin{figure}[ht]
\centering
\[
 [n]=\underbrace{C_\star}_{\text{optional type--1 cycle}}
 \sqcup
 \underbrace{C_1\sqcup\cdots\sqcup C_s}_{\text{type--2 cycles}}
\]
\setlength{\unitlength}{0.72pt}
\begin{picture}(520,145)
  \linethickness{0.8pt}
  \put(20,112){\circle{18}}
  \put(20,112){\makebox(0,0){\scriptsize $i_1$}}
  \put(135,112){\circle{18}}
  \put(135,112){\makebox(0,0){\scriptsize $i_2$}}
  \put(135,64){\circle{18}}
  \put(135,64){\makebox(0,0){\scriptsize $i_3$}}
  \put(20,64){\circle{18}}
  \put(20,64){\makebox(0,0){\scriptsize $i_{r_\star}$}}
  \put(30,112){\vector(1,0){94}}
  \put(135,102){\vector(0,-1){28}}
  \put(124,64){\vector(-1,0){33}}
  \put(77.5,64){\makebox(0,0){$\cdots$}}
  \put(64,64){\vector(-1,0){33}}
  \put(20,74){\vector(0,1){28}}
  \put(77,88){\makebox(0,0){\small $C_\star$}}
  \put(77,42){\makebox(0,0){\small type 1}}
  \put(190,112){\circle{18}}
  \put(190,112){\makebox(0,0){\scriptsize $j_1$}}
  \put(270,112){\circle{18}}
  \put(270,112){\makebox(0,0){\scriptsize $j_2$}}
  \put(270,64){\circle{18}}
  \put(270,64){\makebox(0,0){\scriptsize $j_3$}}
  \put(190,64){\circle{18}}
  \put(190,64){\makebox(0,0){\scriptsize $j_{r_1}$}}
  \put(200,112){\vector(1,0){59}}
  \put(270,102){\vector(0,-1){28}}
  \put(259,64){\vector(-1,0){19}}
  \put(230,64){\makebox(0,0){$\cdots$}}
  \put(220,64){\vector(-1,0){19}}
  \put(190,74){\vector(0,1){28}}
  \put(230,88){\makebox(0,0){\small $C_1$}}
  \put(230,42){\makebox(0,0){\small type 2}}
  \put(330,87){\makebox(0,0){$\cdots$}}
  \put(385,112){\circle{18}}
  \put(385,112){\makebox(0,0){\scriptsize $k_1$}}
  \put(465,112){\circle{18}}
  \put(465,112){\makebox(0,0){\scriptsize $k_2$}}
  \put(465,64){\circle{18}}
  \put(465,64){\makebox(0,0){\scriptsize $k_3$}}
  \put(385,64){\circle{18}}
  \put(385,64){\makebox(0,0){\scriptsize $k_{r_s}$}}
  \put(395,112){\vector(1,0){59}}
  \put(465,102){\vector(0,-1){28}}
  \put(454,64){\vector(-1,0){19}}
  \put(425,64){\makebox(0,0){$\cdots$}}
  \put(415,64){\vector(-1,0){19}}
  \put(385,74){\vector(0,1){28}}
  \put(425,88){\makebox(0,0){\small $C_s$}}
  \put(425,42){\makebox(0,0){\small type 2}}
\end{picture}
\caption{A general summand}
\label{fig:general-partition-weight}
\end{figure}

For $C_\star=(i_1\cdots i_{r_\star})$, each cycle is associated with weight
\[
\begin{aligned}
 W_1(C_\star)
 &=t_M(j_0^{a_{i_1}}\cdots
   j_0^{a_{i_{r_\star}}};q)\Psi_{|C_\star|}(C_\star),\\
 W_2(C_j)
 &=\kk\kappa(a_{C_j})\Phi_{|C_j|}(C_j),\qquad 1\leq j\leq s,
 \end{aligned}
\]
respectively. And total weight is
\[
 W(\alpha)
 =\begin{cases}
\displaystyle
W_1(C_\star)\prod_{j=1}^sW_2(C_j),
& C_\star\text{ present},\\[3ex]
\displaystyle
t_M(1;q)\prod_{j=1}^sW_2(C_j),
& C_\star\text{ absent}.
\end{cases}
\]

\begin{definition}[Graph sum $n$-point Correlators]
\label{def:graph-sums}
The graph sums are
\begin{equation}
 G_{n,M}^{a_1,\ldots,a_n}(z_1,\ldots,z_n;\tau):=
 \sum_{\alpha\in\cP_n^\star}W(\alpha),
\end{equation}
For $n=0$, set $G_{0,M}=t_M(1;q)$.
\end{definition}

For example, the two-point graph sum is
\begin{equation}
 G_2^{ab}(z_1,z_2;\tau)
 =t_M(j_0^a j_0^b;q)\Psi_2(12)
 +t_M(1;q)\kk\kappa(a,b)\Phi_2(12).
 \label{eq:two-point-graph}
\end{equation}
Note that on each finite-dimensional $\g$-module, $\Tr(j_0^a)=0$ because $\g=[\g,\g]$.  It makes the recursive formulas uniform without creating a nonzero one-point graph.

The singular behavior of the graph sum defined here is the same as that predicted by the operator product expansion.
\begin{proposition}
The singular part of graph sum $G_n$ along the $\Delta_{k\ell}=\{(z_1,\ldots,z_n)\in\CC^n| z_k=z_\ell \}$ is given by
\begin{equation}
 \frac{\kk\kappa(a_k,a_\ell)}{(z_k-z_\ell)^2}G_{n-2}^{\widehat k\widehat\ell}
 +\frac{f^{a_ka_\ell}{}_{b}}{z_k-z_\ell}
  G_{n-1}^{b;\widehat k\widehat\ell}.
\label{eq:candidate-principal-part}
\end{equation}
Here $G_{n-2}^{\widehat k\widehat\ell}$ is the $n-2$ points graph sum with $k$-th and $\ell$-th insertions removed. $G_{n-1}^{b;\widehat k\widehat\ell}$ is the $n-1$ points graph sum with $k$-th insertion removed and $\ell$-th insertion replaced by $b$.
\end{proposition}
\begin{proof}

Fix two vertices $k$ and $\ell$.  A graph summand can have a pole at $z_k=z_\ell$ only if $k$ and $\ell$ are consecutive vertices of the same cycle.  If they belong to different cycles, no scalar factor involves their difference; if they belong to the same cycle but are nonadjacent, the definition of Kronecker functions shows that such summands are regular at the diagonal.

Suppose $k$ is immediately followed by $\ell$ in an oriented cycle $C$.  Denote by $C/(k\ell)$ the cycle obtained by collapsing the ordered block $(k,\ell)$ to a fused vertex. Let $w=z_k-z_\ell$, then from  \eqref{eq:H-collision}
\begin{equation*}
 \Res_{w=0}H_r(C;\lambda)
 =\lambda H_{r-1}(C/(k\ell);\lambda).
\end{equation*}
Substituting this identity into the two coefficients give, for $r\geq3$,
\begin{align}
 \Res_{w=0}\Psi_r(C)
 &=\Psi_{r-1}(C/(k\ell)),
 \label{eq:Psi-cycle-residue}\\
 \Res_{w=0}\Phi_r(C)
 &=\Phi_{r-1}(C/(k\ell)).\notag
\end{align}

If the order of $k$ and $\ell$ is reversed, both residues acquire a minus sign:
\begin{align}
 \Res_{w=0}\Psi_r(\cdots \ell,k\cdots)
 &=-\Psi_{r-1},
 \label{eq:Psi-reverse-residue}\\
 \Res_{w=0}\Phi_r(\cdots \ell,k\cdots)
 &=-\Phi_{r-1}.\notag
\end{align}
\paragraph{Simple pole.}
It follows from the local expansions of Kronecker function and definition of $\Phi_n$ and $\Psi_n$ that the simple pole along $\Delta_{k\ell}$ comes from the cycle with length at least $3$ and in which $k$ and $\ell$ are adjacent.

Start with such cycle of type-1
\begin{equation*}
  C_+=(k,\ell,i_3,\ldots,i_r)_1,
  \qquad
  C_-=(\ell,k,i_3,\ldots,i_r)_1.
\end{equation*}

The cycle obtained by collapsing the $k,\ell$ to a fused vertex is denoted by
\begin{equation*}
  \overline C=(b,i_3,\ldots,i_r)_1.
\end{equation*}

Using \eqref{eq:Psi-cycle-residue} and
\eqref{eq:Psi-reverse-residue}, the sum of their residues is
\begin{align*}
&\Res_{w=0}\bigl(W_1(C_+)+W_1(C_-)\bigr)\notag\\
&=\bigl(
 t_M(j_0^{a_k}j_0^{a_\ell}j_0^{a_{i_3}}\cdots j_0^{a_{i_r}};q)
 -t_M(j_0^{a_\ell}j_0^{a_k}j_0^{a_{i_3}}\cdots j_0^{a_{i_r}};q)
\bigr)\Psi_{r-1}(\overline C)\notag\\
&=t_M\!\left(
[j_0^{a_k},j_0^{a_\ell}]
j_0^{a_{i_3}}\cdots j_0^{a_{i_r}};q
\right)\Psi_{r-1}(\overline C)\notag\\
&=f^{a_ka_\ell}{}_{b}
t_M\!\left(
j_0^bj_0^{a_{i_3}}\cdots j_0^{a_{i_r}};q
\right)\Psi_{r-1}(\overline C)\notag\\
&=f^{a_ka_\ell}{}_{b}W_1(\overline C).
\end{align*}

Here $f^{a_ka_\ell}{}_{b}$ are Lie algebra structure constants of $\g$. For the same two cycles $C_+$ and $C_-$ of type-2. The residues are
\begin{align*}
&\Res_{w=0}\bigl(W_2(C_+)+W_2(C_-)\bigr)\notag\\
&=\kk\Bigl(
 \kappa(a_k,a_\ell,a_{i_3},\ldots,a_{i_r})
 -\kappa(a_\ell,a_k,a_{i_3},\ldots,a_{i_r})
 \Bigr)\Phi_{r-1}(\overline C)\notag\\
&=\kk\kappa([a_k,a_\ell],a_{i_3},\ldots,a_{i_r})
  \Phi_{r-1}(\overline C)\notag\\
&=f^{a_ka_\ell}{}_{b}\,
  \kk\kappa(b,a_{i_3},\ldots,a_{i_r})
  \Phi_{r-1}(\overline C)\notag\\
&=f^{a_ka_\ell}{}_{b}W_2(\overline C).
\end{align*}

Let $\cP_n^\star(k\ell)$ denote the subset of $\cP_n^\star$ consists of elements containing cycles of the form $(k,\ell,i_3,\ldots,i_r)$ of either type--1 or type--2. Then for every $\alpha_+\in\cP_n^\star(k\ell)$, there is a unique $\alpha_-\in\cP_n^\star(\ell k)$ obtained by exchanging $k$ and $\ell$. Denote by $\overline \alpha$ the partition by collapsing the $k,\ell$ to a fused vertex. Hence the graph weight contributes to simple pole are
\[
\sum_{\alpha_+\in\cP_n^\star(k\ell)} \Res_{w=0}\bigl(W(\alpha_+)+W(\alpha_-)\bigr)=\sum_{\alpha_+\in\cP_n^\star(k\ell)}f^{a_ka_\ell}{}_{b}W(\overline \alpha).
\]

Fusing the pair $(k\ell)$ gives a bijection between $\cP_n^\star(k\ell)$ and partitions of the remaining $n-1$ labels $\cP_{n-1}^\star$. Hence the coefficient of the order $1$ pole of $G_n$ at $\Delta_{k\ell}$ is 
\[
f^{a_ka_\ell}{}_{b}
  G_{n-1}^{b;\widehat k\widehat\ell}
\]
Here $G_{n-1}^{b;\widehat k\widehat\ell}$ denotes the $n-1$ point function with the $\ell$-th insertion replaced by $b$ and $k$-th insertion removed.

\paragraph{Double pole.}
If $k$ and $\ell$ lie in different cycles, no edge joining them occurs, so the corresponding weight is regular.  If they lie in one cycle of length at least three, at most one edge is singular, so the pole is at most simple. It follows from the local expansions of $\Phi_r$ and $\Psi_r$ that the double pole is only at $r=2$. Example \eqref{eq:explicit-two-kernels} shows that $\Phi_2$ has a double pole while $\Psi_2=(2\pi\ii)^2$ is regular.  Therefore a double pole can arise only when $(k\ell)$ is an independent
type--2 two-cycle.

Its weight is
\begin{equation*}
  \kk\kappa(a_k,a_\ell)\Phi_2(k\ell)
  =\kk\kappa(a_k,a_\ell)
   \left(\frac1{(z_k-z_\ell)^2}+O(1)\right).
\end{equation*}
Analogously to the discussion for type--1. Deleting the pair $(k\ell)$ gives a bijection between partitions containing this pair and partitions of the remaining $n-2$ labels. Hence the coefficient of the order $2$ of $G_n$ pole at $\Delta_{k\ell}$ is
\begin{equation*}
 \kk\kappa(a_k,a_\ell)G_{n-2}^{\widehat k\widehat\ell}.
\end{equation*}
Here $G_{n-2}^{\widehat k\widehat\ell}$ denotes the $n-2$ point function with the $k$-th and $\ell$-th insertions removed.

 Summing over all lower partitions shows the singular part of $G_n$ along $\Delta_{k\ell}$ is
\begin{equation}
\frac{\kk\kappa(a_k,a_\ell)}{(z_k-z_\ell)^2}G_{n-2}^{\widehat k\widehat\ell}
 +\frac{f^{a_ka_\ell}{}_{b}}{z_k-z_\ell}
  G_{n-1}^{b;\widehat k\widehat\ell}
\end{equation}
as expected.
\end{proof}

\section{The Graph-Sum $n$-Correlators Equal the 
	$n$-Point Trace}
\label{sec:formal-proof}

In this section we prove the correlation functions of affine currents are exactly given by our graph sum. We show that $S_n=G_n$ as a formal power series in $\CC[[u_1^{\pm1},\ldots,u_n^{\pm1};q]]$. We first define two formal power series $T^S$ and $T^G$ on arbitrary affine mode words. They are the coefficients in the expansion in $u_i$ of the trace and the graph sum respectively. We then show that they satisfy the same four characterizing properties which will finally show they are equal.

\subsection{Formal traces}

A mode word is an expression
\begin{equation*}
 w=j_{m_1}^{b_1}\cdots j_{m_r}^{b_r},
 \qquad b_i\in\g,\quad m_i\in\ZZ,
\end{equation*}
and its total mode number is $m(w):=\sum_i m_i$.  Define
\begin{equation}
 T^S_M(w;q):=
 \begin{cases}
 \displaystyle
 \sum_{N\geq0}\Tr_{M_N}(w)q^N,&m(w)=0,\\[1.2ex]
 0,&m(w)\neq0.
 \end{cases}
 \label{eq:TS-definition}
\end{equation}
If $m(w)=0$, the operator $w$ preserves every homogeneous subspace, so
each coefficient is an ordinary finite-dimensional trace.  For a word
$w_0=j_0^{b_1}\cdots j_0^{b_r}$ consisting only of zero modes, write
\begin{equation}
 t_M(j_0^{b_1}\cdots j_0^{b_r};q):=
 \sum_{N\geq0}\Tr_{M_N}(w_0)q^N.
 \label{eq:t-zero-word}
\end{equation}
Thus $T^S_M(w_0;q)=t_M(j_0^{b_1}\cdots j_0^{b_r};q)$.  We usually omit the subscript $M$.

These $T^S(w;q)$ serve as coefficients of expansion of correlation functions in $u_i$
\begin{equation}
 S_n^{a_1\cdots a_n}(u_1,\ldots,u_n;q)
 =(2\pi\ii)^n
 \sum_{\substack{m_i\in\ZZ\\m_1+\cdots+m_n=0}}
 T^S(j_{m_1}^{a_1}\cdots j_{m_n}^{a_n};q)
 \prod_{i=1}^n u_i^{-m_i}.
 \label{eq:S-formal-expansion}
\end{equation}
It is an element of
$\CC[[u_1^{\pm1},\ldots,u_n^{\pm1};q]]$.

We collect some basic properties of $T^S$ which will be used later.
\begin{proposition}
\label{prop:TS-relations}
The family $T^S_M(w;q)$, defined for all finite mode words, has the
following properties.
\begin{enumerate}[label=\textup{(\roman*)},leftmargin=2.7em]
 \item It vanishes if $m(w)\neq0$.
 \item Adjacent factors satisfy
 \begin{align}
 &T^S(\cdots j_k^a j_\ell^b\cdots;q)
 -T^S(\cdots j_\ell^b j_k^a\cdots;q)\notag\\
 &\quad=T^S(\cdots j_{k+\ell}^{[a,b]}\cdots;q)
 +\kk k\kappa(a,b)\delta_{k+\ell,0}
 T^S(\cdots\widehat{j_k^a}\widehat{j_\ell^b}\cdots;q).
 \label{eq:TS-exchange}
 \end{align}
 \item If $w=j_m^a j_{m_2}^{b_2}\cdots j_{m_r}^{b_r}$ has total mode number zero, then
 \begin{equation}
  T^S(j_{m_2}^{b_2}\cdots j_{m_r}^{b_r}j_m^a;q)
  =q^mT^S(j_m^a j_{m_2}^{b_2}\cdots j_{m_r}^{b_r};q)
  \qquad\text{in }\CC((q)).
  \label{eq:TS-cyclicity}
 \end{equation}
 \item On zero-mode words,
 \begin{equation}
  T^S(j_0^{b_1}\cdots j_0^{b_r};q)
  =t_M(j_0^{b_1}\cdots j_0^{b_r};q).
  \label{eq:TS-zero-data}
 \end{equation}
\end{enumerate}
\end{proposition}

\begin{proof}
The first and fourth assertions are definitions, and the second is the
operator identity of modes Lie algebra \eqref{eq:affine-bracket-universal}.  For (iii), $j_{m_2}^{b_2}\cdots j_{m_r}^{b_r}$ has total mode number $-m$.
The maps
\begin{equation*}
 j_{m_2}^{b_2}\cdots j_{m_r}^{b_r}:M_N\longrightarrow M_{N+m},\qquad
 j_m^a:M_{N+m}\longrightarrow M_N
\end{equation*}
are between finite-dimensional spaces.  The trace identity gives
\begin{equation*}
 \Tr_{M_N}(j_m^a j_{m_2}^{b_2}\cdots j_{m_r}^{b_r})
 =\Tr_{M_{N+m}}(j_{m_2}^{b_2}\cdots j_{m_r}^{b_r}j_m^a).
\end{equation*}
Multiply by $q^{N+m}$ and \eqref{eq:TS-cyclicity} follows.
\end{proof}

\subsection{Formal graph sum}

In this section, we collect some properties of graph sum and its expansion coefficients which will be used to prove trace and graph sum have the same expansion.

\begin{lemma}
\label{lem:kernel-identities}
Let $C=(i_1\cdots i_r)$ be a cycle of length $r$.  The type-1 elliptic cycle coefficient $\Psi_r$ and
the type-2 elliptic cycle coefficient $\Phi_r$ have the following properties.
\begin{enumerate}[label=\textup{(\roman*)},leftmargin=2.7em]

\item They are elliptic in each $u_i$, equivalently invariant under $u_i\mapsto qu_i$.
 
 \item They are invariant under simultaneous scaling of all multiplicative coordinates, equivalently
 \[
 \Phi(cu_{i_1},\ldots,cu_{i_r})= \Phi(u_{i_1},\ldots,u_{i_r}), \quad \Psi(cu_{i_1},\ldots,cu_{i_r})= \Psi(u_{i_1},\ldots,u_{i_r}),\quad \forall c\in\CC^*.
 \] 
 
 \item Their expansions in the domain $|u_{i_1}|>|u_{i_2}|>\cdots >|u_{i_r}|$ have no negative powers
 of $q$ so define elements of
 \begin{equation*}
  \CC[[u_{i_1}^{\pm1},\ldots,u_{i_r}^{\pm1};q]].
 \end{equation*}

 \item For $C=(i_1\cdots i_r)$, write
 $[u_{i_1}^0\cdots u_{i_r}^0]$ for extraction of the constant term
 in these variables from the radial expansion in the domain $  0 < |q| < |u_{i_n}| < \dots < |u_{i_1}| < 1 $.
 More precisely, if $j_1<\cdots<j_r$ are the labels of $C$ in increasing
 order, then, for a test function $F(u_{i_1},\ldots,u_{i_r})$, this means
 \[
 [u_{i_1}^0\cdots u_{i_r}^0]\,
 \iota_{|u_{j_1}|>\cdots>|u_{j_r}|}F.
 \]
 This radial order depends on the ambient variable labels, rather than
 on the cyclic order of $C$. In the coefficient extractions below,
 this radial expansion is understood. Then
 \begin{equation}
 \begin{aligned}
  [u_{i_1}^0\cdots u_{i_r}^0]\Phi_r(C)&=0,\\
  [u_{i_1}^0\cdots u_{i_r}^0]\Psi_r(C)&=
  \begin{cases}
   (2\pi\ii)^r,& i_1<i_2<\cdots<i_r \text{ up to cyclic rotation},\\
   0,&\text{otherwise}.
  \end{cases}
 \end{aligned}
  \label{eq:kernel-constant-terms}
 \end{equation}
\end{enumerate}
\end{lemma}

\begin{proof}
The quasiperiodicity \eqref{eq:Kronecker-periods} cancels between the incoming and outgoing edge at each vertex, proving ellipticity.  This establishes (i).

The proof of (ii) follows immediately from the definition.

Put $u=\e^{2\pi\ii z}$ and $v=\e^{2\pi\ii\lambda}$.  The Jacobi triple product in our setting \eqref{eq:odd-theta}
\[
\theta(z,\tau)=iq^{\frac18}u^{-\frac{1}{2}}(1-u)\prod_{s\geq1}(1-q^s)(1-uq^s)(1-u^{-1}q^s)
\] 
gives
\begin{equation}
 K(z,\lambda;\tau)
 =(2\pi\ii)\frac{uv-1}{(u-1)(v-1)} \prod_{s\geq1} \frac{(1-q^s)^2(1-uvq^s)(1-u^{-1}v^{-1}q^s)} {(1-uq^s)(1-u^{-1}q^s)(1-vq^s)(1-v^{-1}q^s)}.
 \label{eq:K-product}
\end{equation}
After geometric expansion of every denominator, a fixed $q$-coefficient receives
contributions from only finitely many factors, and no negative power of $q$ occurs. The coefficient extractions in \eqref{eq:Psi-definition} and \eqref{eq:Phi-definition} preserve nonnegative $q$-degree.  Expanding the remaining rational functions in the domain $  0 < |q| < |u_n| < \dots < |u_1| < 1 $ gives an element of $\CC[[u_1^{\pm1},\ldots,u_n^{\pm1};q]]$. 

We give the constant-term calculation in full.  The two radial Fourier expansions of the Kronecker function are
\begin{align}
 K(z,\lambda)
 &=-2\pi\ii\sum_{m\in\ZZ}\frac{u^m}{1-vq^m},
 &&|q|<|u|<1,
 \label{eq:K-Fourier-inner}\\
 K(z,\lambda)
 &=-2\pi\ii\sum_{m\in\ZZ}
 \frac{vq^m u^m}{1-vq^m},
 &&1<|u|<|q|^{-1}.
 \label{eq:K-Fourier-outer}
\end{align}
For $C=(i_1\cdots i_r)$, let
\begin{equation*}
 d(C):=\#\{s:i_s<i_{s+1}\},\qquad i_{r+1}=i_1.
\end{equation*}
Thus $1\leq d(C)\leq r-1$.  In the cyclic product, all vertex powers
cancel only when all edge Fourier indices are equal.  With
$D=2\pi\ii\lambda$, this yields
\begin{equation}
 [u_{i_1}^0\cdots u_{i_r}^0] H_r(C;\lambda)
 =(-D)^r\sum_{m\in\ZZ}
 \frac{(vq^m)^{d(C)}}{(1-vq^m)^r}.
 \label{eq:CT-H}
\end{equation}
For $m\neq0$, the summand is regular at $D=0$ and is multiplied by
$D^r$, so it cannot contribute to the coefficient of $D^{r-1}$ defining
$\Psi_r$.  The $m=0$ term gives
\begin{align}
 [u_{i_1}^0\cdots u_{i_r}^0]\Psi_r(C)
 &=(2\pi\ii)^r[D^{r-1}]
 \left(\frac{D}{\e^D-1}\right)^r
 \e^{(d(C)+1)D}\notag\\
 &=\frac{(2\pi\ii)^r}{(r-1)!}
 B_{r-1}^{(r)}(d(C)+1).
 \label{eq:Psi-CT-Bernoulli}
\end{align}
By \cref{prop:bernoulli-factorization},
\begin{equation*}
 B_{r-1}^{(r)}(x)=\prod_{j=1}^{r-1}(x-j).
\end{equation*}
The value in \eqref{eq:Psi-CT-Bernoulli} is therefore nonzero precisely
when $d(C)=r-1$, which is the increasing cyclic order, and then equals
$(2\pi\ii)^r$.

It remains to treat $\Phi_r$.  Set
\begin{equation*}
 \Sigma_{r,d}(D):=\sum_{m\in\ZZ}
 \frac{(vq^m)^d}{(1-vq^m)^r},
 \qquad 1\leq d\leq r-1.
\end{equation*}
The rational functions
\begin{equation*}
 E_s(t):=(t\partial_t)^{s-1}\frac{t}{(1-t)^2},
 \qquad 1\leq s\leq r-1,
\end{equation*}
form a basis of
\begin{equation*}
 \left\{\frac{t p(t)}{(1-t)^r}:\deg p\leq r-2\right\}.
\end{equation*}
Indeed, this space has dimension $r-1$, and the functions $E_s$ are
linearly independent because their pole orders at $t=1$ are the distinct
integers $s+1$, for $1\leq s\leq r-1$.  Hence there are constants
$c_{r,d,s}$ such that
\begin{equation*}
 \frac{t^d}{(1-t)^r}=\sum_{s=1}^{r-1}c_{r,d,s}E_s(t).
\end{equation*}
Using the Fourier expansion of $P_2$ gives
\begin{equation}
 \Sigma_{r,d}(D)=\sum_{s=1}^{r-1}
 \frac{c_{r,d,s}}{(2\pi\ii)^{s+1}}P_2^{(s-1)}(\lambda).
 \label{eq:Sigma-P2}
\end{equation}
For every $s\geq0$,
\begin{equation}
 \Res_{\lambda=0}
 (P_1(\lambda)+\pi\ii)P_2^{(s)}(\lambda)\,\dd\lambda=0.
 \label{eq:residue-zero-P1P2}
\end{equation}
For odd $s$ the integrand is even.  For $s=2\ell$, its only two possible
residue contributions are
\begin{equation*}
 -(2\ell+1)!G_{2\ell+2}
 +(2\ell+1)!G_{2\ell+2}=0,
\end{equation*}
by the expansions following \eqref{eq:P1P2}.  Finally,
$\theta'/\theta=-(P_1+\pi\ii)$, so
\begin{align*}
 [u_{i_1}^0\cdots u_{i_r}^0]\Phi_r(C)
 &=(-2\pi\ii)^r\Res_{\lambda=0}
 (P_1(\lambda)+\pi\ii)
 \Sigma_{r,d(C)}(D)\,\dd\lambda=0
\end{align*}
by \eqref{eq:Sigma-P2} and \eqref{eq:residue-zero-P1P2}.  This completes the proof of (iv).
\end{proof}

\subsection{The graph coefficient \texorpdfstring{$T^G$}{TG}}

 Let $w=j_{m_1}^{b_1}\cdots j_{m_r}^{b_r}$ be an arbitrary mode word, and let $G_r^{b_1,\ldots,b_r}(u_1,\ldots, u_r;q)$ denote the graph
sum of \cref{def:graph-sums}. 

Define the graph coefficient by taking the coefficient in the radial
expansion determined by the displayed order of the word:
\begin{equation}
 T^G_M(w;q):=\frac1{(2\pi\ii)^r}
 \left[u_1^{-m_1}u_2^{-m_2}\cdots u_r^{-m_r}\right]
 \iota_{\scrD_r}G_r^{b_1,\ldots,b_r}(u_1,u_2,\ldots,u_r;q)
 \label{eq:TG-definition}
\end{equation}
where  $\scrD_r$ is the domain $  0 < |q| < |u_r| < \dots < |u_1| < 1 $. Thus the first factor of the word is assigned to the outermost variable, the second factor to the next variable, and so on.  If the factors of the word are permuted, both the labels and their coordinate slots are permuted,
and the coefficient is then taken in the radial order of the new word.

Since $\Phi$, $\Psi$, and the formal zero-mode traces have no negative
powers of $q$, every fixed $q$-coefficient of $G_r$ is a finite sum of
products of finite-dimensional trace coefficients and Laurent
coefficients of meromorphic functions of the $u_i$.
Lemma~\ref{lem:kernel-identities}(i) therefore guarantees that
$T^G_M(w;q)\in\CC[[q]]$.

\begin{proposition}
\label{prop:TG-relations}
The family $T^G_M(w;q)$ satisfies the four assertions of \cref{prop:TS-relations}.
\end{proposition}

\begin{proof}
The graph sum is invariant under simultaneous permutation of
label-coordinate pairs.  Explicitly, for every $\sigma\in S_r$,
\begin{equation*}
 G_r^{b_{\sigma(1)},\ldots,b_{\sigma(r)}}
 (u_{\sigma(1)},\ldots,u_{\sigma(r)};q)
 =G_r^{b_1,\ldots,b_r}(u_1,\ldots,u_r;q).
\end{equation*}
This follows easily from the symmetry of the graph sum.

\smallskip
\noindent\emph{Property (i): total mode number.}
Both $\Phi$ and $\Psi$ depend only on differences $z_i-z_j$, and hence are
unchanged when all multiplicative coordinates are simultaneously scaled:
\begin{equation*}
 G_r(cu_1,\ldots,cu_r;q)=G_r(u_1,\ldots,u_r;q).
\end{equation*}
If a monomial $u_1^{-m_1}\cdots u_r^{-m_r}$ occurs in the ordered
expansion, simultaneous scaling multiplies it by
$c^{-\sum_i m_i}$.  Uniqueness of Laurent coefficients therefore forces
its coefficient to vanish unless $\sum_i m_i=0$.  By
\eqref{eq:TG-definition}, this proves
$T^G_M(w;q)=0$ whenever $m(w)\neq0$.

\smallskip
\noindent\emph{Property (ii): adjacent exchange.}
By definition, 
\begin{align*}
 &T^G(\cdots j_k^a j_\ell^b\cdots;q)
 -T^G(\cdots j_\ell^b j_k^a\cdots;q)\\
 = &\frac{1}{(2\pi\ii)^r}[\cdots u_i^{-k}u_{i+1}^{-\ell}\cdots]\iota_{\cdots >|u_i|>|u_{i+1}|>\cdots}G^{\cdots a,b\cdots}(\ldots, u_i,u_{i+1},\ldots)\\
 & -\frac{1}{(2\pi\ii)^r}[\cdots u_i^{-\ell}u_{i+1}^{-k}\cdots]\iota_{\cdots >|u_i|>|u_{i+1}|>\cdots}G^{\cdots b,a\cdots}(\ldots, u_i,u_{i+1},\ldots)\\
  =&  \frac{1}{(2\pi\ii)^r}[\cdots u_i^{-k}u_{i+1}^{-\ell}\cdots]\left(\iota_{|u_i|>|u_{i+1}|}-\iota_{|u_{i+1}|>|u_{i}|}\right)G^{\cdots a,b\cdots}(\ldots, u_i,u_{i+1},\ldots).
\end{align*}
Here $j^a_k$ is at the $i$-th position of this word. In the last equality, we simultaneously exchange insertion and position, then interchange the dummy indices $u_i$ and $u_{i+1}$. We know that only singular terms in $u_i/u_{i+1}$ will have different expansions and therefore contribute to the above expression.

Recall we have proved that the singular part of $G_r$ along $\Delta_{i,i+1}$ is
\begin{equation}
\frac{\kk\kappa(a,b)}{(z_i-z_{i+1})^2}G^{\widehat{a}\widehat{b}}_{r-2}+\frac1{z_i-z_{i+1}}
 \left(G^{[a,b];\widehat{a}\widehat{b}}_{r-1}\right).
 \label{eq:graph-additive-OPE}
\end{equation}
In multiplicative coordinates $u=\e^{2\pi\ii z}$, the singular part is
\begin{equation}
\frac{(2\pi\ii)^2\kk\kappa(a,b)u_i u_{i+1}}{(u_i-u_{i+1})^2}
 G_{r-2} +\frac{(2\pi\ii)u_i}{u_i-u_{i+1}}G_{r-1}.
 \label{eq:graph-multiplicative-OPE}
\end{equation}

We now compare the two ordered expansions. Set
\begin{equation*}
 \delta(v/u):=\sum_{s\in\ZZ}(v/u)^s.
\end{equation*}
The elementary formal-distribution identities
\begin{align*}
 \iota_{|u|>|v|}\frac{u}{u-v}
 -\iota_{|v|>|u|}\frac{u}{u-v}
 &=\delta(v/u),\\
 \iota_{|u|>|v|}\frac{uv}{(u-v)^2}
 -\iota_{|v|>|u|}\frac{uv}{(u-v)^2}
 &=(v\partial_v)\delta(v/u)
\end{align*}
follow by expanding the first expression in nonnegative powers of $v/u$
and the second in nonnegative powers of $u/v$.  The regular part of the graph sum has the same expansion on
the two sides of this one diagonal.  Hence \eqref{eq:graph-multiplicative-OPE} gives that
\begin{align*}
& \left(\iota_{|u_i|>|u_{i+1}|}-\iota_{|u_{i+1}|>|u_{i}|}\right)G^{\cdots a,b\cdots}(\ldots, u_i,u_{i+1},\ldots) \\
= & (2\pi\ii)\delta \left(\frac{u_{i+1}}{u_i}\right)G^{\cdots \widehat{a},[a,b]\cdots}_{r-1}(\ldots,\widehat{u_i},u_{i+1},\ldots)
\\ &\quad+(2\pi\ii)^2\kk\kappa(a,b) \left( u_{i+1}\frac{\partial}{\partial u_{i+1}}\right)\delta\left(\frac{u_{i+1}}{u_i}\right)G^{\cdots,\widehat{a},\widehat{b}\cdots}_{r-2}(\ldots,\widehat{u_i},\widehat{u_{i+1}},\ldots).
\end{align*}

Taking the coefficient of $u^{-k}v^{-\ell}$ uses
\begin{align*}
 [u^{-k}v^{-\ell}]\,\delta(v/u)F(v)
 &=[v^{-(k+\ell)}]F(v),\\
 [u^{-k}v^{-\ell}]\,(v\partial_v)\delta(v/u)
 &=k\,\delta_{k+\ell,0}.
\end{align*}
Applying these two identities, we obtain
\begin{align*}
 &T^G(\cdots j_k^a j_\ell^b\cdots;q)
 -T^G(\cdots j_\ell^b j_k^a\cdots;q)\\
 &\quad=T^G(\cdots j_{k+\ell}^{[a,b]}\cdots;q)
 +\kk k\kappa(a,b)\delta_{k+\ell,0}
 T^G(\cdots\widehat{j_k^a}\widehat{j_\ell^b}\cdots;q),
\end{align*}
which is property (ii).

\smallskip
\noindent\emph{Property (iii): graded cyclicity.}
Let $w=j_m^a j_{m_2}^{b_2}\cdots j_{m_r}^{b_r}$ have total mode number zero, and write the coordinates of $j_{m_2}^{b_2}\cdots j_{m_r}^{b_r}$ as $u_2,\ldots,u_r$. 
For $r\geq 2$, the coefficient for $j_m^a j_{m_2}^{b_2}\cdots j_{m_r}^{b_r}$ can be taken in the relative radial chamber
\[
 |u|>|u_2|>\cdots>|u_r|>|qu|.
\]
After permuting the insertion-coordinate pairs, the coefficient for $j_{m_2}^{b_2}\cdots j_{m_r}^{b_r}j_m^a$ is taken from the same graph sum in the chamber
\[
 |u_2|>\cdots>|u_r|>|u|>|qu_2|.
\]
These chambers contain the respective ordered chambers of \eqref{eq:TG-definition}. The coordinate-dependent factors are holomorphic there, so their Laurent expansions restrict to the same expansions in the smaller chambers. The case $r=1$ has $m=0$ and is immediate.

Substituting $u\mapsto qu$ sends the first relative chamber onto the second. By ellipticity in the $u$-coordinate,
\begin{equation*}
 G_r^{a,b_2,\ldots,b_r}(qu,u_2,\ldots,u_r;q)
 =G_r^{a,b_2,\ldots,b_r}(u,u_2,\ldots,u_r;q).
\end{equation*}
Taking the coefficient of $u^{-m}u_2^{-m_2}\cdots u_r^{-m_r}$ on the left introduces a factor $q^{-m}$ from the substitution. Consequently,
\begin{equation*}
 T^G(j_{m_2}^{b_2}\cdots j_{m_r}^{b_r}j_m^a;q)
 =q^mT^G(j_m^a j_{m_2}^{b_2}\cdots j_{m_r}^{b_r};q).
\end{equation*}
When $m<0$, this identity is read in $\CC((q))$, exactly as in
\cref{prop:TS-relations}; both $T^G$ coefficients themselves still lie in
$\CC[[q]]$.

\smallskip
\noindent\emph{Property (iv): zero-mode data.}
Set $m_1=\cdots=m_r=0$.  Definition \eqref{eq:TG-definition} then asks
for the constant term of each graph summand.  Since different cycles use disjoint sets of variables, this constant term is the product of the constant terms of the individual type-1 and type-2 elliptic cycle coefficients $\Psi$ and $\Phi$.  By \eqref{eq:kernel-constant-terms}, a type--2 cycle contributes zero. Therefore every partition containing a type--2 cycle vanishes, including every partition having no type--1 cycle.

The only possible surviving partition consists of a single type--1 cycle containing all vertices.  Among its cyclic orders, the only survivor is the increasing order $(1\,2\cdots r)$, whose weight is exactly $t_M(j_0^{b_1}\cdots j_0^{b_r};q)$. In summary,
\begin{equation*}
 T^G(j_0^{b_1}\cdots j_0^{b_r};q)
 =t_M(j_0^{b_1}\cdots j_0^{b_r};q).
\end{equation*}
For the empty word, the same statement is the convention
$T^G(1;q)=t_M(1;q)$.  This proves property (iv) and completes the proof of all four assertions.
\end{proof}

\subsection{Proof of the main theorem}
Now we can state and prove our main result of this paper.

\begin{theorem}
\label{thm:formal-cycle}
For every module with finite dimensional homogeneous subspace and every $n\geq0$,
\begin{equation}
S_{n,M}=G_{n,M}\quad\text{in }
 \CC[[u_1^{\pm1},\ldots,u_n^{\pm1};q]]
 \label{eq:formal-cycle-equality}
\end{equation}
Equivalently, for $n\geq1$,
\begin{align}
 &(2\pi\ii)^n
 \sum_{\substack{\mbf\in\ZZ^n\\\sum_i m_i=0}}
 T^S_M(j_{m_1}^{a_1}\cdots j_{m_n}^{a_n};q)
 \prod_{i=1}^n u_i^{-m_i}\notag\\
 &\qquad=\iota_{\scrD_n}
 \sum_{\alpha\in\cP_n^\star}
 \begin{cases}
 \displaystyle
 t_M(j_0^{a_{i_1}}\cdots j_0^{a_{i_r}};q)
 \Psi_r(C_\star)
 \prod_{C\neq C_\star}\kk\kappa(a_C)\Phi_{|C|}(C),
 &C_\star=(i_1\cdots i_r)\text{ present},\\[3ex]
 \displaystyle
 t_M(1;q)\prod_{C\in\alpha}\kk\kappa(a_C)\Phi_{|C|}(C),
 &C_\star\text{ absent}.
 \end{cases}
 \label{eq:formal-cycle-expanded}
\end{align}
\end{theorem}

\begin{proof}
By \cref{prop:TS-relations,prop:TG-relations}, $T^S_M$ and $T^G_M$
satisfy the same properties. Now we only need to prove they agree on every mode word.

Put $D=T^S-T^G$.  We prove
$D(w;q)=0$ by induction on the length $r$ of $w$.  The value of $D$ on empty word is zero because its value is prescribed.  Assume the assertion at all lengths smaller than $r$.  We have proved adjacent-exchange can be reduced to two terms with length $r-1$ and $r-2$. Both vanish for $D$ by the induction hypothesis.  Therefore, at length $r$,
\begin{equation}
 D(\cdots j_k^a j_\ell^b\cdots;q)
 =D(\cdots j_\ell^b j_k^a\cdots;q);
 \label{eq:D-permutation}
\end{equation}
the difference is invariant under every permutation of its mode factors.

If all modes are zero, the prescribed initial condition gives $D=0$.
Otherwise, use \eqref{eq:D-permutation} to put a nonzero factor
$j_m^a$ at the front and write $w=j_m^a j_{m_2}^{b_2}\cdots j_{m_r}^{b_r}$.  Moving it to the end by
adjacent exchanges leaves $D$ unchanged, while graded cyclicity gives
\begin{equation*}
 D(j_{m_2}^{b_2}\cdots j_{m_r}^{b_r}j_m^a;q)
 =q^mD(j_m^a j_{m_2}^{b_2}\cdots j_{m_r}^{b_r};q).
\end{equation*}
Consequently $(1-q^m)D(w;q)=0$ in $\CC((q))$.  Since $m\neq0$, the
factor $1-q^m$ is invertible, so $D(w;q)=0$.  This closes the induction.
\end{proof}

\section{Convergence and Meromorphic Continuation}
\label{sec:convergence}

Let $ \widetilde {M}_{\kk} $ and $M$ be the  $\hat \g $-modules as in Section 1.
In this section we prove the convergence and meromorphic continuation of $n$-point trace.

\subsection{Convergence of zero-mode traces}

\begin{lemma}
	$	\Tr_M q^{d } =   \sum_{ N = 0 }^{\infty} {\rm dim} \,  M_N q^N $ converges on $ |q | < 1 $.
\end{lemma}

\begin{proof}
	Since $  {\rm dim} \, M_N \leq   {\rm dim} \,  (\widetilde{M}_{\kk})_N $, it is enough to prove the result for the case $ M = \widetilde M_{\kk }  $.
	By PBW Theorem, we have 
	\[ 	\Tr_{\widetilde M_{\kk}} q^{d } =   {\rm dim } \, M_{\rm top}    \, \prod_{ j =1 }^\infty  ( 1 -q^j )^{-{\rm dim } \, \g } \]
	which is convergent.
\end{proof}

\begin{lemma}
	\label{lem:finite-dimensional-trace-bound} Let $E$ be a finite dimensional $\g $-module,
	let $ U(\g)_{\leq r}$ denote the filtration of universal enveloping algebra $ U(\g)$ which is spanned by products of at most $r$ elements of $\g$. For every $w\in U(\g)_{\leq r}$, there is a constant $C_w$, independent of the $E$, such that
	\begin{equation}
	\left|\Tr_E\rho_E(w)\right|
	\leq C_w(\dim E)^{r+1}.
	\label{eq:finite-dimensional-trace-bound}
	\end{equation}
\end{lemma}

\begin{proof}
	Write $d=\dim E$. Fix a compact real form of the simple Lie algebra $\mathfrak g$ and equip $E$ with an invariant positive-definite Hermitian form. Choose a Chevalley basis
	\[
	\mathcal B
	=
	\{h_i\}_{i=1}^{\operatorname{rank}\mathfrak g}
	\cup
	\{e_\alpha\}_{\alpha\in\Delta}
	\]
	of $\mathfrak g$, here $\Delta$ is the root system. Every element of $\mathcal B$ belongs to a root $\mathfrak{sl}_2$-subalgebra.  More precisely, for every root $\alpha$, the elements
	\[
	e_\alpha,\qquad
	e_{-\alpha},\qquad
	h_\alpha=[e_\alpha,e_{-\alpha}]
	\]
	span a copy of $\mathfrak{sl}_2$.  Restricting $E$ to this
	$\mathfrak{sl}_2$, we obtain an orthogonal decomposition into
	irreducible $\mathfrak{sl}_2$-modules.
	
	Recall that the irreducible $\mathfrak{sl}_2$-module $L(m)$ has
	dimension $m+1$ and admits an orthonormal basis
	$v_0,\ldots,v_m$ such that
	\[
	hv_j=(m-2j)v_j,
	\]
	\[
	ev_j=\sqrt{j(m-j+1)}\,v_{j-1},
	\qquad
	fv_j=\sqrt{(j+1)(m-j)}\,v_{j+1}.
	\]
	It follows that
	\[
	\|\rho_{L(m)}(h)\|_{\mathrm{op}}\leq m,
	\qquad
	\|\rho_{L(m)}(e)\|_{\mathrm{op}},
	\ \|\rho_{L(m)}(f)\|_{\mathrm{op}}
	\leq \frac{m+1}{2}.
	\]
	Every irreducible $\mathfrak{sl}_2$-summand of $E$ has dimension at most $d$.  Therefore, after fixing Chevalley basis, there exists a constant $C_{\mathfrak g}$, depending only on $\mathfrak g$, such that
	\[
	\|\rho_E(b)\|_{\mathrm{op}}
	\leq
	C_{\mathfrak g}d
	\]
	for every $b\in\mathcal B$ and every finite-dimensional
	$\mathfrak g$-module $E$.
	
	By the definition of the PBW filtration, $w$ can be written as
	\[
	w
	=
	\sum_{\nu=1}^{N}
	c_\nu
	b_{\nu,1}\cdots b_{\nu,s_\nu},
	\qquad
	b_{\nu,j}\in\mathcal B,
	\qquad
	0\leq s_\nu\leq r.
	\]
	Using the triangle inequality and the submultiplicativity of the
	operator norm, we obtain
	\[
	\begin{aligned}
	\|\rho_E(w)\|_{\mathrm{op}}
	&\leq
	\sum_{\nu=1}^{N}
	|c_\nu|
	\prod_{j=1}^{s_\nu}
	\|\rho_E(b_{\nu,j})\|_{\mathrm{op}}\\
	&\leq
	\sum_{\nu=1}^{N}
	|c_\nu|
	(C_{\mathfrak g}d)^{s_\nu}\\
	&\leq
	C_w' d^r,
	\end{aligned}
	\]
	where
	\[
	C_w'
	=
	\sum_{\nu=1}^{N}
	|c_\nu|C_{\mathfrak g}^{s_\nu}
	\]
	depends only on $w$ and not on $E$.
	
	Finally, for every operator $A$ on a $d$-dimensional Hermitian
	space,
	\[
	|\Tr A|
	\leq
	d\|A\|_{\mathrm{op}}.
	\]
	Hence
	\[
	\left|
	\Tr_E\rho_E(w)
	\right|
	\leq
	d\|\rho_E(w)\|_{\mathrm{op}}
	\leq
	C_w'd^{r+1}.
	\]
	Taking $C_w=C_w'$ proves the claim.  The case $E=0$ is immediate.
\end{proof}

We will use the following equivalent description of the convergence of series to show that the convergence of the character controls the convergence of the zero mode trace.
\begin{proposition}[Zero-mode convergence]
	\label{prop:zero-mode-convergence}
	For $ a_1 , \dots , a_r \in \g $, the trace 
	$t_M(j_0^{a_1}\cdots j_0^{a_r};q)$  
	converges absolutely and locally uniformly in the domain $\{|q|<1\}$.
\end{proposition}
\begin{proof}
	Since the radius of the convergence of the series $ \sum_{ N= 0 } {\rm dim} \, M_N q^N $ is greater or equal to $1$, we have
	$ \overline {\lim} ({\rm dim} \, M_N )^{\frac 1 N} \leq  1$, this implies that for every $\epsilon > 0 $,
	$ {\rm dim} \, M_N \leq C_\epsilon e^{\epsilon N } $ for some constant $C_\epsilon $. So we have
	\[ \sum_{ N=0}^\infty | \Tr_{M_N}  j^{a_1}_0  \cdots j^{a_r}_0 | \, |q|^N \leq \sum_{ N =0 }^\infty 
	C \, ({\rm dim}\, M_N )^{r+1}   |q|^N \leq \sum_{N=0}^\infty  CC_\epsilon^{r+1} e^{( r+1 ) \epsilon N} |q|^N \]
	This converges on $ | q | < e^{ - ( r+1 ) \epsilon} $, since $\epsilon$ is arbitrary, the results follow. 
\end{proof}

\subsection{Convergence of current traces}
Now we rewrite the graph formula as a finite linear combination
\begin{equation}
 G_{n,M}(\ubf;q)
 =\sum_{\alpha=1}^{A_n}
 f_{n,\alpha}(\ubf;q)t_M(j_0^{b_{\alpha,1}}\cdots j_0^{b_{\alpha,r_\alpha}};q),
 \qquad b_{\alpha,j}\in\g,\quad 0\leq r_\alpha\leq n.
 \label{eq:finite-zero-mode-reduction}
\end{equation}
Each $f_{n,\alpha}$ is a finite product of $\Phi$ and at most one $\Psi$ together with level $\kk$ and extended Killing form. Since $f_{n,\alpha}$ is already meromorphic, the convergence of current correlation is totally controlled by zero mode traces.

\begin{theorem}[Convergence of current traces]
\label{thm:current-convergence}
For every $n\geq0$ and $a_1,\ldots,a_n\in\g$, the formal power series of trace \eqref{eq:S-formal-expansion} converges absolutely and locally uniformly on $\scrD_n$ to graph sum $G_{n,M}$, and the same finite graph expression gives the unique meromorphic continuation, hence defines a holomorphic function on
\begin{equation*}
 \operatorname{Conf}_n(E_\tau)
 :=\{(z_1,\ldots,z_n)\in E_\tau^n:z_i\neq z_j\text{ for }i\neq j\}.
\end{equation*}

\end{theorem}

\begin{proof}
By \cref{prop:zero-mode-convergence}, every series in the finite sum \eqref{eq:finite-zero-mode-reduction} converges normally on $|q|<1$.  The triple-product expression \eqref{eq:K-product} shows that the scalar kernels are jointly holomorphic at $q=0$ in the ordered configuration away from collisions.  Thus
$G_{n,M}$ is holomorphic on the Reinhardt domain 
\begin{equation*}
 \scrD_n^0:=
 \{(q,\ubf):|q|<|u_n|<\cdots<|u_1|<1\}.
\end{equation*}
For fixed $q$, its Laurent coefficients are given by iterated Cauchy integrals on any product of circles $|u_i|=\rho_i$ satisfying $|q|<\rho_n<\cdots<\rho_1<1$.  Normal convergence permits termwise multiplication and integration. By \cref{thm:formal-cycle}, these coefficients are exactly the formal current-trace coefficients in \eqref{eq:S-formal-expansion}.

The Laurent expansion theorem on logarithmically convex Reinhardt
domains now implies absolute and locally uniform convergence of this
Laurent series to $G_{n,M}$ on $\scrD_n$; see \cite[Section~1.5]{Scheidemann}.  

The scalar kernels are meromorphic and elliptic in every insertion, with
poles only when $z_i-z_j\in\ZZ+\tau\ZZ$.  Hence the finite graph formula
continues meromorphically to the elliptic configuration space. 

\end{proof}

\subsection{General insertion}
 In this section, we discuss the convergence and meromorphic continuation problem of general $n-$point functions. We will first consider the one point case and show it controls any point case.

Let $M=\bigoplus_{N\geq0}M_N$ be a graded module.  For homogeneous $v_1,\ldots,v_n$, define the formal trace
\begin{equation}
 S_{n,M}^{v_1,\ldots,v_n}(\ubf;q)
 :=\left(\prod_{i=1}^n(2\pi\ii u_i)^{|v_i|}\right)
 \Tr_M\!\left(Y_M(v_1,u_1)\cdots Y_M(v_n,u_n)q^d\right),
 \label{eq:general-insertion-trace}
\end{equation}
For current states $v_i=j_{-1}^{a_i}\one$, it agrees with the current trace defined earlier.

\begin{lemma}[Convergence of one-point traces]
\label{lem:arbitrary-affine-one-point}
For every $v\in V^{\kk}(\g)$, let $|v|$ denote its affine degree and set $o(v)=v_{(|v|-1)}$, the series
\begin{equation}
 S_M(v;q)
 =\sum_{N\geq0}
 \Tr_{M_N}\!\left(o(v)\right)q^N
 \label{eq:arbitrary-affine-one-point}
\end{equation}
converges absolutely and locally uniformly for $|q|<1$. 
\end{lemma}

\begin{proof}
It is enough to consider a homogeneous element
\begin{equation}
 v=j_{-\ell_1}^{a_1}\cdots j_{-\ell_r}^{a_r}\one,
 \qquad
 \ell_\alpha\geq1, \quad
 |v|=\ell_1+\cdots+\ell_r.
 \label{eq:zhu-section-pbw-word}
\end{equation}
The case $r=0$ is exactly the character.

We claim that its one-point trace is the following iterated residues of the $r$-point graph sum.
\begin{align}
 S_M(v;q)
 = & \frac{u^{|v|}}{(2\pi i)^r}
 \Res_{u_1=u}\cdots\Res_{u_r=u}
 \left[
  \prod_{i=1}^{r}
  (u_i-u)^{-\ell_i}
   G_{r,M}^{a_1,\ldots,a_r}
  (u_1,\ldots,u_r;q)
 \prod_{i=1}^{r}\frac{\dd u_i}{u_i}\right].
 \label{eq:one-point-as-current-coefficient}
\end{align}
Here $G_{r,M}$ is the graph sum and the residues are taken in the domain
\begin{equation*}
 |u_1-u|>|u_2-u|>\cdots> |u_r-u|>0.
\end{equation*}

We prove \eqref{eq:one-point-as-current-coefficient} by induction on the word length $r$.  For $r=0$, the state is $\one$ and both sides are the character $t_M(1;q)$.  Suppose the formula is known for a word $v'$ of length $r-1$ and write
\begin{equation*}
 v=j_{-\ell}^{a}v'.
\end{equation*}
The operator product expansion says that $Y_M(j_{-\ell}^{a}v',u')$ is the local coefficient of $j^a(u)Y_M(v',u')$ selected by $(u-u')^{-\ell}$.  This is exactly \eqref{eq:one-point-as-current-coefficient} for $v$. 

Convergence now follows easily.  By \eqref{eq:finite-zero-mode-reduction},
\begin{equation*}
 G_{r,M}(\ubf;q)
 =\sum_{\text{finite}}
 f_{r,\alpha}(\ubf;q)t_M(j_0^{b_{\alpha,1}}\cdots j_0^{b_{\alpha,r_\alpha}};q),
 \qquad b_{\alpha,j}\in\g,\quad 0\leq r_\alpha\leq r,
\end{equation*}
is a finite sum.  Every zero-mode series $t_M(j_0^{b_{\alpha,1}}\cdots j_0^{b_{\alpha,r_\alpha}};q)$ converges absolutely and locally uniformly for $|q|<1$ by \cref{prop:zero-mode-convergence}; the elliptic cycle coefficients are normally convergent on compact subsets avoiding elliptic diagonals. Taking finite iterated residue map preserves normal convergence.  Hence the right-hand side of \eqref{eq:one-point-as-current-coefficient} converges absolutely and locally uniformly for $|q|<1$.
\end{proof}

The general case is given by the recurrent formula \eqref{prop:zhu-recursion-recalled}, which reduces $n-$point functions to fewer insertion functions. 

We first introduce square bracket formalism. Given a vertex operator algebra $(V,Y(\cdot,z),\one,\omega)$, we define a new state-field correspondence map
\begin{equation*}
Y[a,z]=Y(a,e^{2\pi\ii z}-1)e^{2\pi \ii z|a|}.
\end{equation*}
This gives a new vertex operator algebra structure on $V$ with identity element $\one$ and Virasoro element $(2\pi\ii)^2(\omega-\frac{c\one}{24})$. If we write $Y[a,z]=\sum_{n\in\ZZ}a_{[n]}z^{-n-1}$, then 
\begin{equation*}
a_{[m]}=(2\pi\ii)^{-m-1}\Res_z \left(Y(a,z)\left(\ln(1+z)\right)^m(1+z)^{|a|-1} \right)
\end{equation*}
With the help of square bracket formalism, we may state the recurrent formula \eqref{prop:zhu-recursion-recalled} as follows.
\begin{proposition}[Proposition 4.4.2 of \cite{Zhu}]
\label{prop:zhu-recursion-recalled}
Let $a,a_1,\ldots,a_n$ be homogeneous vectors of a vertex operator
algebra and let $M$ be a lower-bounded module with finite-dimensional
homogeneous subspaces. We have \cite{Zhu}
\begin{align}
 &S_{M,n+1}
 ((a,w),(a_1,z_1),\ldots,(a_n,z_n);\tau)
 \notag\\
 &={}
S_{M,n}
 ((a[-1]a_1,z_1),(a_2,z_2),\ldots,(a_n,z_n);\tau)
 \notag\\
 &\quad-
 \sum_{k\geq2}G_{2k}(\tau)\,
 S_{M,n}
 ((a[2k-1]a_1,z_1),(a_2,z_2),\ldots,(a_n,z_n);\tau)
 \notag\\
 &\quad+
 \sum_{m\geq0}\wp_{m+1}(w-z_1,\tau)\,
 S_{M,n}
 ((a[m]a_1,z_1),(a_2,z_2),\ldots,(a_n,z_n);\tau)
 \notag\\
 &\quad+
 \sum_{j=2}^{n}\sum_{m\geq0}
 \bigl(
   \wp_{m+1}(w-z_j,\tau)
   -\wp_{m+1}(z_1-z_j,\tau)
 \bigr)
 \notag\\[-0.2em]
 &\times
 S_{M,n}
 ((a_1,z_1),\ldots,(a[m]a_j,z_j),\ldots,
  (a_n,z_n);\tau).
 \label{eq:zhu-additive-recursion}
\end{align}
For fixed vectors $b,c$, the truncation property gives $b[m]c=0$ for all
sufficiently large $m$. So every sum in \eqref{eq:zhu-additive-recursion} is finite.
\end{proposition}

\begin{remark}
Although Zhu's recursion is formulated in the presence of a Virasoro element, its derivation uses only the vertex-algebra identities, the grading, and graded cyclicity of the trace. Hence it remains valid as a formal identity for the critical-level affine vertex algebra $V^{-h^\vee}(\g)$. Here traces are taken against $q^d$, and the square-bracket fields are defined by $Y[a,z]=Y(a,e^{2\pi\ii z}-1)e^{2\pi\ii z|a|}$ for homogeneous $a$.
\end{remark}

The essential feature of \eqref{eq:zhu-additive-recursion} is that every trace on the right has one fewer insertion than the trace on the left. Iterating the formula therefore gives the following immediate consequence.

\begin{corollary}
\label{cor:zhu-reduction-one-point}
For fixed $v_1,\ldots,v_n$, the formal $n$-point trace is a finite sum
\begin{equation}
 S_{M,n}
 ((v_1,z_1),\ldots,(v_n,z_n);\tau)
 =\sum_{\text{finite}}
 C_\alpha(\zbf;\tau)\,S_M(w_\alpha;q),
 \label{eq:zhu-finite-one-point-reduction}
\end{equation}
where $S_M(w;q)$ is the one point function. Each $w_\alpha$ is obtained from the $v_i$ by finitely many square-bracket products, and each $C_\alpha$ is a finite sum of finite products of $G_{2k}(\tau)$ and $\wp_m(z_i-z_j,\tau)$.  Only finitely many states $w_\alpha$ occur.
\end{corollary}

\begin{theorem}[Convergence of arbitrary correlation]
\label{thm:zhu-recursion-general-convergence}
Let
\begin{equation*}
 M=\bigoplus_{N=0}^{\infty}M_N
\end{equation*}
be a lower-bounded $V^{\kk}(\g)$-module with finite-dimensional homogeneous subspaces, and assume that its character converges in the domain $\{0<|q|<1 \}$.  Then for every $n\geq1$ and every
$v_1,\ldots,v_n\in V^{\kk}(\g)$, the trace
\begin{equation}
 S_M^{v_1,\cdots,v_n}(u_1;\ldots;u_n;q)
 \label{eq:zhu-section-general-trace}
\end{equation}
converges absolutely and locally uniformly on
\begin{equation*}
1>|u_1|>\cdots>|u_n|>|q|>0.
\end{equation*}
to an analytic function and has a unique meromorphic continuation to $\operatorname{Conf}_n(E_\tau)$.  
\end{theorem}

\begin{proof}
 The case $n=1$ is just \cref{lem:arbitrary-affine-one-point}.

The $n$-point trace is a combination of the convergent one-point series supplied by \cref{lem:arbitrary-affine-one-point}, with normally convergent scalar coefficients.  This makes explicit that, after one-point convergence is known, no independent multi-point estimate remains.

Replacing the ordered Fourier expansions of the scalar coefficients by the meromorphic functions $G_{2k}(\tau)$ and $\wp_m(z_i-z_j,\tau)$ gives a meromorphic continuation away from the elliptic diagonals.  Zhu's recurrent formula, together with graded cyclicity, identifies the expressions obtained from neighboring lifts, so they descend to a meromorphic function on
$\operatorname{Conf}_n(E_\tau)$.  It agrees with \eqref{eq:zhu-section-general-trace} on the initial chamber and is unique
by the identity theorem.
\end{proof}

\section{Low-Point Formulas}
\label{sec:low-points}

This section evaluates the type-1 and type-2 elliptic cycle coefficients and the correlation functions explicitly for $n\leq 4$.

\subsection{Cycle coefficients}
By convention,
\begin{equation*}
 \Phi_1=0,\qquad \Psi_1=2\pi\ii.
\end{equation*}
For two points, put $z=z_1-z_2$.  The sigma addition formula gives
\begin{equation*}
 \frac{\sigma(\lambda+z)\sigma(\lambda-z)}
 {\sigma(\lambda)^2\sigma(z)\sigma(-z)}
 =\wp_2(\lambda)-\wp_2(z).
\end{equation*}
Therefore
\begin{equation*}
 H_2(z,-z;\lambda)
 =\lambda^2\bigl(\wp_2(\lambda)-\wp_2(z)\bigr)
 =1-\wp_2(z)\lambda^2+3G_4\lambda^4+O(\lambda^6),
\end{equation*}
and
\begin{equation}
\Psi_2=(2\pi\ii)^2,\qquad
 \Phi_2(12)=\wp_2(z_1-z_2)+G_2.
 \label{eq:explicit-two-kernels}
\end{equation}

For $C=(123)$, set
\begin{align*}
 \mu_1&=z_1-z_2,&
 \mu_2&=z_2-z_3,&
 \mu_3&=z_3-z_1,\\
 \mathcal Z_{123}&:=\sum_{s=1}^3\wp_1(\mu_s),&
 \mathcal W_{123}&:=\sum_{s=1}^3\wp_2(\mu_s),&
 \mathcal W'_{123}&:=\sum_{s=1}^3\wp_2'(\mu_s).
\end{align*}
Equation \eqref{eq:H-Bell} gives
\begin{align*}
 H_{3,0}&=1,\qquad
 H_{3,1}=\mathcal Z_{123},\\
 H_{3,2}&=\frac12
 \bigl(\mathcal Z_{123}^2-\mathcal W_{123}\bigr)=0,\\
 H_{3,3}&=-\frac13\mathcal Z_{123}^3
 -\frac16\mathcal W'_{123}.
\end{align*}
The vanishing of $H_{3,2}$ is \eqref{eq:H-vanishing}, equivalently the
Frobenius--Stickelberger identity
$\mathcal Z_{123}^2=\mathcal W_{123}$.  Hence
\begin{align}
 \Psi_3(123)
 &=\frac{(2\pi\ii)^3}{2}
 +(2\pi\ii)^2\mathcal Z_{123},
 \label{eq:Psi3}\\
 \Phi_3(123)
 &=G_2\mathcal Z_{123}
 +\frac13\mathcal Z_{123}^3
 +\frac16\mathcal W'_{123}.
 \label{eq:Phi3}
\end{align}
The formulas for any other directed three-cycle follow by permutation of
the variables.

For $C=(1234)$, put
\begin{equation*}
 (\mu_1,\mu_2,\mu_3,\mu_4)
 =(z_1-z_2,z_2-z_3,z_3-z_4,z_4-z_1)
\end{equation*}
and abbreviate
\begin{align*}
 \mathcal Z&:=\sum_{s=1}^4\wp_1(\mu_s),&
 \mathcal W&:=\sum_{s=1}^4\wp_2(\mu_s),\\
 \mathcal W_1&:=\sum_{s=1}^4\wp_2'(\mu_s),&
 \mathcal W_2&:=\sum_{s=1}^4\wp_2''(\mu_s).
\end{align*}
Then
\begin{align*}
 H_{4,1}&=\mathcal Z,\\
 H_{4,2}&=\frac12(\mathcal Z^2-\mathcal W),\\
 H_{4,3}&=\frac16
 (\mathcal Z^3-3\mathcal Z\mathcal W-\mathcal W_1)=0,\\
 H_{4,4}&=\frac1{24}\left(
 -3\mathcal Z^4+6\mathcal Z^2\mathcal W
 +3\mathcal W^2-\mathcal W_2\right)+G_4.
\end{align*}
In the last line we used
$\mathcal W_1=\mathcal Z^3-3\mathcal Z\mathcal W$, which is the preceding
vanishing identity.  Substitution into \eqref{eq:Psi-H} and
\eqref{eq:Phi-H} gives
\begin{align}
 \Psi_4(1234)
 &=\frac{(2\pi\ii)^4}{6}
 +\frac{(2\pi\ii)^3}{2}\mathcal Z
 +\frac{(2\pi\ii)^2}{2}(\mathcal Z^2-\mathcal W),
 \label{eq:Psi4}\\
 \Phi_4(1234)
 &=\frac{
 3\mathcal Z^4-6\mathcal Z^2\mathcal W
 -3\mathcal W^2+\mathcal W_2}{24}
 +\frac{G_2}{2}(\mathcal Z^2-\mathcal W).
 \label{eq:Phi4}
\end{align}

\subsection{Trace functions}

The zero- and one-point traces are
\begin{equation*}
 S_0=t_M(1;q),\qquad S_1^a=0.
\end{equation*}

The one-point vanishing follows grade by grade because every finite
dimensional $\g$-module has zero trace on $\g$.

The two-point formula is
\begin{equation}
 S_2^{ab}(z_1,z_2;\tau)
 =(2\pi\ii)^2t_M(j_0^aj_0^b;q)
 +t_M(1;q)\kk\kappa(a,b)
 \bigl(\wp_2(z_1-z_2)+G_2\bigr).
 \label{eq:S2-explicit}
\end{equation}

For three points, abbreviate
$\Psi_3(ijk)=\Psi_3(z_i,z_j,z_k)$ and similarly for $\Phi_3$.  The full
nonzero formula is
\begin{align}
 S_3^{a_1a_2a_3}
={}&t_M(j_0^{a_1}j_0^{a_2}j_0^{a_3};q)\Psi_3(123)
 +t_M(j_0^{a_1}j_0^{a_3}j_0^{a_2};q)\Psi_3(132)\notag\\
 &+t_M(1;q)\kk\kappa(a_1,a_2,a_3)\Phi_3(123)
 +t_M(1;q)\kk\kappa(a_1,a_3,a_2)\Phi_3(132).
 \label{eq:S3-explicit}
\end{align}

Let
\begin{equation*}
 \mathcal C_4
 :=\{(1234),(1243),(1324),(1342),(1423),(1432)\}
\end{equation*}
be the six directed four-cycles modulo cyclic rotation, and set
\begin{equation*}
 \mathcal P_{2+2}:=\{(12|34),(13|24),(14|23)\}.
\end{equation*}
Then
\begin{align}
 S_4^{a_1a_2a_3a_4}
={}&\sum_{\rho\in\mathcal C_4}
 \Bigl[
 t_M(j_0^{a_{\rho_1}}j_0^{a_{\rho_2}}
 j_0^{a_{\rho_3}}j_0^{a_{\rho_4}};q)\Psi_4(\rho)
 \notag\\
 &\hspace{3.2cm}
 +t_M(1;q)\kk\kappa(a_{\rho_1},a_{\rho_2},a_{\rho_3},a_{\rho_4})
 \Phi_4(\rho)\Bigr]\notag\\
 &+\sum_{(ij|k\ell)\in\mathcal P_{2+2}}
 \Bigl[
 t_M(1;q)\kk\kappa(a_i,a_j)\kk\kappa(a_k,a_\ell)
 \Phi_2(ij)\Phi_2(k\ell)\notag\\
 &\hspace{2.0cm}
 +t_M(j_0^{a_i}j_0^{a_j};q)\Psi_2(ij)
 \kk\kappa(a_k,a_\ell)\Phi_2(k\ell)\notag\\
 &\hspace{2.0cm}
 +t_M(j_0^{a_k}j_0^{a_\ell};q)\Psi_2(k\ell)
 \kk\kappa(a_i,a_j)\Phi_2(ij)\Bigr].
 \label{eq:S4-explicit}
\end{align}

\section{Nodal Degeneration}
\label{sec:nodal-degeneration}

This section studies the limit $\tau\to \ii\infty$.  The elliptic curve degenerates to a nodal rational curve, the cycle coefficients $\Phi_n$ and $\Psi_n$ tend to $\Phi_n^{(0)}$ and $\Psi_n^{(0)}$, and the zero mode trace reduces to the ground-state space.

\subsection{Limits of the cycle coefficients}
The two presentations
\begin{equation*}
 E_\tau=\CC/(\ZZ+\tau\ZZ)\cong\CC^*/q^\ZZ,
 \qquad u=\e^{2\pi\ii z},
\end{equation*}
identify the translations by $1$ and $\tau$ with $u\mapsto u$ and
$u\mapsto qu$.  As $\tau\to\ii\infty$, or equivalently $q\to 0$, the quotient becomes a rational curve with one node.  Its normalization is
\begin{equation*}
 \nu:\mathbb P^1_u\longrightarrow E_0,
\end{equation*}
and $u=0,\infty$ are the two points above the node.  The ordered radial domain becomes
\begin{equation}
 1>|u_1|>\cdots>|u_n|>0,
 \label{eq:plane-radial-chamber}
\end{equation}
the standard radial order on the sphere.

The Jacobi triple product applied to \eqref{eq:odd-theta} gives
\begin{equation}
 \theta(z,\tau)=2q^{1/8}\sin(\pi z)
 \prod_{m\geq1}(1-q^m)(1-uq^m)(1-u^{-1}q^m).
 \label{eq:theta-product-cusp}
\end{equation}
Therefore
\begin{equation}
 \theta(z,\tau)=2q^{1/8}\sin(\pi z)(1+O(q)),
 \qquad
 \theta'(0,\tau)=2\pi q^{1/8}(1+O(q)).
 \label{eq:theta-cusp}
\end{equation}
Also,
\begin{equation*}
 G_{2r}(\tau)
 =2\zeta(2r)
 +\frac{2(2\pi\ii)^{2r}}{(2r-1)!}
 \sum_{m\geq1}\frac{m^{2r-1}q^m}{1-q^m},
\end{equation*}
where $\zeta$ is the Riemann zeta function.  In particular,
$G_2\to\pi^2/3$.

It follows from \eqref{eq:theta-cusp} that
\begin{equation}
 K(z,\lambda;\tau)=K^{(0)}(z,\lambda)+O(q),
 \qquad
 K^{(0)}(z,\lambda)
 =\pi\frac{\sin\pi(z+\lambda)}
 {\sin\pi z\,\sin\pi\lambda}
 =\pi(\cot\pi\lambda+\cot\pi z).
 \label{eq:K0}
\end{equation}
For later comparison,
\begin{equation}
 \wp_2(z,\tau)+G_2(\tau)\longrightarrow
 \pi^2\csc^2(\pi z)
 =(2\pi\ii)^2\frac{u_i u_j}{(u_i-u_j)^2},
 \qquad u=\frac{u_i}{u_j}.
 \label{eq:P2-cusp}
\end{equation}

For $C=(i_1\cdots i_n)$, put
\begin{equation*}
 v_s:=\e^{2\pi\ii(z_{i_s}-z_{i_{s+1}})}
 =\frac{u_{i_s}}{u_{i_{s+1}}},
 \qquad i_{n+1}=i_1,
\end{equation*}
and let $D=2\pi\ii\lambda$.  Since
\begin{equation*}
 \cot(\pi z)=\ii\frac{\e^{2\pi\ii z}+1}
 {\e^{2\pi\ii z}-1},
\end{equation*}
equation \eqref{eq:K0} becomes
\begin{equation}
 K^{(0)}(z_{i_s}-z_{i_{s+1}},\lambda)
 =2\pi\ii\frac{v_s\e^D-1}{(v_s-1)(\e^D-1)}.
 \label{eq:K0-multiplicative}
\end{equation}
Define
\begin{equation}
 H_n^{(0)}(C;\lambda)
 :=\lambda^n\prod_{s=1}^n
 K^{(0)}(z_{i_s}-z_{i_{s+1}},\lambda)
 =D^n\prod_{s=1}^n
 \frac{v_s\e^D-1}{(v_s-1)(\e^D-1)}.
 \label{eq:H0}
\end{equation}
Coefficient extraction commutes with the $q\to0$ limit on compact sets,
so define
\begin{align}
 \Psi_n^{(0)}(C)&:=\lim_{q\to0}\Psi_n(C;\tau),&
 \Phi_n^{(0)}(C)&:=\lim_{q\to0}\Phi_n(C;\tau).
 \label{eq:kernel-zero-limits}
\end{align}

\begin{proposition}[Cusp kernels]
\label{prop:cusp-kernels}
For every directed cycle $C=(i_1\cdots i_n)$ of length $n\geq2$,
\begin{align}
 \Psi_n^{(0)}(C)
 &=(2\pi\ii)^n\sum_{t=1}^n
 \prod_{\substack{1\leq s\leq n\\s\neq t}}
 \frac{u_{i_s}}{u_{i_s}-u_{i_{s+1}}},
 \label{eq:Psi0}\\
 \Phi_n^{(0)}(C)
 &=-(2\pi\ii)^n\prod_{s=1}^n
 \frac{u_{i_s}}{u_{i_s}-u_{i_{s+1}}}.
 \label{eq:Phi0}
\end{align}
The conventions $\Psi_1^{(0)}=2\pi\ii$ and $\Phi_1^{(0)}=0$ remain in
force.
\end{proposition}

\begin{proof}
By \eqref{eq:Psi-definition} and \eqref{eq:H0},
\begin{align*}
 \Psi_n^{(0)}(C)
 &=(2\pi\ii)^n\Res_{D=0}
 \e^D\prod_{s=1}^n
 \frac{v_s\e^D-1}{(v_s-1)(\e^D-1)}\,\dd D.
\end{align*}
Set $x=\e^D$.  Since $\e^D\dd D=\dd x$, this becomes
\begin{equation*}
 (2\pi\ii)^n\Res_{x=1}
 \frac{\prod_s(v_sx-1)}
 {(x-1)^n\prod_s(v_s-1)}\,\dd x.
\end{equation*}
With $y=x-1$,
\begin{equation*}
 \frac{v_sx-1}{v_s-1}
 =1+\frac{u_{i_s}}{u_{i_s}-u_{i_{s+1}}}y.
\end{equation*}
The residue selects degree $n-1$ in the product, giving
\eqref{eq:Psi0}.

Next,
\begin{equation*}
 \lambda\frac{\theta'(\lambda;\tau)}{\theta(\lambda;\tau)}
 \longrightarrow
 \pi\lambda\cot(\pi\lambda)
 =\frac D2\frac{\e^D+1}{\e^D-1}.
\end{equation*}
The definition of $\Phi_n$ therefore gives
\begin{equation*}
 \Phi_n^{(0)}(C)
 =-\frac{(2\pi\ii)^n}{2}
 [y^n]\left(
 \frac{y+2}{y+1}
 \prod_{s=1}^n
 \left(1+\frac{u_{i_s}}{u_{i_s}-u_{i_{s+1}}}y\right)\right).
\end{equation*}
Write the product as $P(y)=\sum_{j=0}^np_jy^j$.  Since
\begin{equation*}
 [y^n]\frac{y+2}{y+1}P(y)
 =p_n+\sum_{j=0}^n(-1)^{n-j}p_j
 =p_n+(-1)^nP(-1),
\end{equation*}
and $\prod_sv_s=1$, one has
\begin{equation*}
 p_n=\prod_{s=1}^n\frac{u_{i_s}}{u_{i_s}-u_{i_{s+1}}},
 \qquad P(-1)=(-1)^np_n.
\end{equation*}
Thus the extracted coefficient is $2p_n$, proving \eqref{eq:Phi0}.
\end{proof}

\subsection{Plane cycle coefficients}

We first recall the genus-zero correlation functions.  For $v\in M_0$, regard $v'\in M_0^*$ as an element of the restricted dual $M'=\bigoplus_{N=0}^{\infty} M^{*}_N$. For $a_1,\ldots,a_n\in\g$, write
\begin{equation}
 F^{a_1,\ldots,a_n}_n(v',v;\ubf)
 :=\left\langle v',
 j^{a_1}(u_1)\cdots j^{a_n}(u_n)v\right\rangle,
 \qquad F_0(v',v):=\langle v',v\rangle.
 \label{eq:plane-matrix-coefficient}
\end{equation}
Here $j^a(u)$ is the plane current in \eqref{eq:plane-current}.
The product is expanded in radial order $|u_1|>\cdots>|u_n|>0$.
Geometrically, $v$ and $v'$ are inserted at $0$ and $\infty$ on
$\mathbb P^1_u$.

Following \cite{FrenkelZhu}, a \textbf{chain} is a nonempty subset of $[n]$ equipped with a linear order.  We denote it by $C'=\langle i_1,\ldots,i_r\rangle$, with $r\geq1$.  A cycle $C=(i_1\cdots i_r)$ has the same meaning as before.

Let $\cP_n^{\text{pl}}$ be the set of partitions of $[n]$ into
disjoint cycles and \textbf{at most one chain}. Thus a partition has the form
\begin{equation*}
 [n]=\underbrace{C'_\star}_{\text{optional chain}}
 \sqcup C_1\sqcup\cdots\sqcup C_s.
\end{equation*}
For $n=0$, the set consists of the empty partition.

Define the chain coefficient and the cycle coefficient analogously by
\begin{align}
 \Psi_r^{(\mathrm{pl})}(C')
 &:=\frac1{(u_{i_1}-u_{i_2})\cdots
 (u_{i_{r-1}}-u_{i_r})u_{i_r}},
 \qquad r\geq1,
 \label{eq:plane-Psi}\\
 \Phi_r^{(\mathrm{pl})}(C)
 &:=-\frac1{(u_{i_1}-u_{i_2})\cdots
 (u_{i_r}-u_{i_1})},
 \qquad r\geq2.
 \label{eq:plane-Phi}
\end{align}
The special case when $r=1$ is $\Psi_1^{(\mathrm{pl})}(\langle i\rangle)=u_i^{-1}$.  We also retain the convention $\Phi_1^{(\mathrm{pl})}=0$ because there is no cycle of length $1$.

The two component weights are
\begin{align*}
 W_1^{(\mathrm{pl})}(C';v',v)
 &:=\left\langle v',
 j_0^{a_{i_1}}\cdots j_0^{a_{i_r}}v\right\rangle
 \Psi_r^{(\mathrm{pl})}(C'),\\
 W_2^{(\mathrm{pl})}(C)
 &:=\kk\kappa(a_{i_1},\ldots,a_{i_r})
 \Phi_r^{(\mathrm{pl})}(C).
\end{align*}
The extended Killing tensor is exactly \eqref{eq:extended-form}, with
the cyclic order of $C$.  For a partition $\alpha$, set
\begin{equation}
 W_{v',v}^{(\mathrm{pl})}(\alpha)
 :=\begin{cases}
 \displaystyle
 W_1^{(\mathrm{pl})}(C'_\star;v',v)
 \prod_{C\in\alpha,\,C\neq C'_\star}
 W_2^{(\mathrm{pl})}(C),
 &C'_\star\text{ present},\\[2ex]
 \displaystyle
 \langle v',v\rangle\prod_{C\in\alpha}
 W_2^{(\mathrm{pl})}(C),
 &C'_\star\text{ absent}.
 \end{cases}
 \label{eq:plane-partition-weight}
\end{equation}
In particular, a partition with no chain still carries the pairing
$\langle v',v\rangle$, the matrix coefficient of the empty zero-mode
word.

\begin{proposition}[Genus-zero correlation functions \cite{FrenkelZhu}]
\label{prop:plane-chain-cycle}
For every $v\in M_0$, $v'\in M_0^*$, and $a_1,\ldots,a_n\in\g$, the formal matrix coefficient \eqref{eq:plane-matrix-coefficient} converges absolutely and locally uniformly in $|u_1|>\cdots>|u_n|>0$ to
\begin{equation}
 F^{a_1,\ldots,a_n}_n(v',v;\ubf)
 =\sum_{\alpha\in\cP_n^{(\mathrm{pl})}}
 W_{v',v}^{(\mathrm{pl})}(\alpha).
 \label{eq:plane-chain-cycle-formula}
\end{equation}
The right-hand side is rational, with poles only at $u_i=0$ and $u_i=u_j$.
\end{proposition}

\subsection{Degeneration relation}

We now relate these plane coefficients to the cycle coefficients defined earlier.  For a cycle $C'=\langle i_1\cdots i_r\rangle$, define cyclic rotation operator $\rho$ by
\begin{equation*}
 \rho(C'):=\langle i_{2},i_{3},\ldots,i_{r},i_1\rangle.
\end{equation*}
Define the cyclic summation of $\Psi_r^{\mathrm{pl}}$ by
\begin{equation}
\Psi_r^{\mathrm{pl,cyc}}=\sum_{t=1}^r\Psi_r^{(\mathrm{pl})}(\rho^t(C')) =\sum_{t=1}^r\frac1{u_{i_t}}
 \prod_{\substack{1\leq s\leq r\\s\neq t}}
 \frac1{u_{i_s}-u_{i_{s+1}}}
\end{equation}

Choose a basis $\{v_\beta\}$ of the finite-dimensional ground-state space $M_0$ and its dual basis $\{v^\beta\}$ in $M_0^*$. Define the trace of the matrix coefficients on ground-state as
\begin{equation}
 S_{n,M_0}^{(\mathrm{pl}),a_1\cdots a_n}(\ubf)
 :=\sum_\beta F^{a_1\cdots a_n}_n(v^\beta,v_\beta;\ubf)
 =\sum_\beta
 \left\langle v^\beta,
 j^{a_1}(u_1)\cdots j^{a_n}(u_n)v_\beta
 \right\rangle
 \label{eq:plane-ground-block}
\end{equation}
in the domain $1>|u_1|>\cdots>|u_n|>0$.  This is independent of the chosen basis. For $n=0$, it gives $S_{0,M_0}^{(\mathrm{pl})}=\dim M_0$.

Summing \eqref{eq:plane-chain-cycle-formula} over the dual bases gives
\begin{equation}
 S_{n,M_0}^{(\mathrm{pl})}(\ubf)
 =\sum_{\alpha\in\cP_n^{(\mathrm{pl})}}
 W_{M_0}^{(\mathrm{pl})}(\alpha),
 \label{eq:plane-ground-cycle-formula}
\end{equation}
The weights are
\begin{equation}
 W_{M_0}^{(\mathrm{pl})}(\alpha)
 :=\begin{cases}
 \displaystyle
 \Tr_{M_0}(j_0^{a_{i_1}}\cdots j_0^{a_{i_r}})
 \Psi_r^{(\mathrm{pl})}(C'_{\star})
 \prod_{C\in\alpha,\,C\neq C'_{\star}}
 \kk\kappa(a_C)\Phi_{|C|}^{(\mathrm{pl})}(C),
 &C'_\star\text{ present},\\[2ex]
 \displaystyle
 \dim M_0\prod_{C\in\alpha}
 \kk\kappa(a_C)\Phi_{|C|}^{(\mathrm{pl})}(C),
 &C'_{\star}\text{ absent}.
 \end{cases}
 \label{eq:plane-ground-partition-weight}
\end{equation}
Here $\kappa(a_C)$ means the same extended Killing form as in the torus graph sum.  Although a trace of zero modes is invariant under cyclic rotation, the plane chain coefficient is not; all linear orders in $\cP_n^{(\mathrm{pl})}$ must still be summed.

Note that for any $C'$, all $\rho^t(C'_{\star}),0\leq t\leq r-1$ give the same $\Tr_{M_0}(j_0^{a_{i_1}}\cdots j_0^{a_{i_r}})$. Sum over all such cyclic rotation, we have
\begin{equation}
 S_{n,M_0}^{(\mathrm{pl})}(\ubf)
 =\sum_{\alpha\in\cP_n^{(\mathrm{pl})}}
 W_{M_0}^{(\mathrm{pl})}(\alpha)=\sum_{\alpha\in\cP_n^{\star}}
 W_{M_0}^{(\mathrm{pl,cyc})}(\alpha)
\end{equation}
with weight
\begin{equation*}
 W_{M_0}^{(\mathrm{pl,cyc})}(\alpha)
 :=\begin{cases}
 \displaystyle
 \Tr_{M_0}(j_0^{a_{i_1}}\cdots j_0^{a_{i_r}})
 \Psi_r^{(\mathrm{pl,cyc})}(C'_{\star})
 \prod_{C\in\alpha,\,C\neq C'_{\star}}
 \kk\kappa(a_C)\Phi_{|C|}^{(\mathrm{pl})}(C),
 &C'_\star\text{ present},\\[2ex]
 \displaystyle
 \dim M_0\prod_{C\in\alpha}
 \kk\kappa(a_C)\Phi_{|C|}^{(\mathrm{pl})}(C),
 &C'_{\star}\text{ absent}.
 \end{cases}
\end{equation*}

\begin{theorem}
\label{thm:cusp-cycle-formula}
Let $M$ be a $V^{\kk}(\g)$-module with finite dimensional homogeneous subspaces. For every $a_1,\ldots,a_n\in\g$,  the $\tau\to i\infty$ limit of correlation functions on the elliptic curves is given by
\begin{equation}
\lim_{\tau\to\ii\infty}S_{n,M}^{a_1\cdots a_n}(\zbf;\tau)
 =\left(\prod_{i=1}^n2\pi\ii\,u_i\right)
 S_{n,M_0}^{(\mathrm{pl}),a_1\cdots a_n}(\ubf).
 \label{eq:cusp-plane-trace-limit}
\end{equation}
Equivalently, in $q$-coordinate,
\begin{equation}
 S_{n,M}(\zbf;\tau)=S_{n,M}^{(0)}(\ubf)+O(q)=\prod_{i=1}^n (2\pi\ii u_i)S_{n,M_0}^{\mathrm{pl}}(\ubf)+O(q)
 \label{eq:normalized-cusp-trace}
\end{equation}
Moreover, $S_{n,M}^{(0)}$ is the limit obtained by the replacements
\begin{equation}
 t_M(j_0^{a_{i_1}}\cdots j_0^{a_{i_r}};q)
 \longmapsto
 \Tr_{M_0}(j_0^{a_{i_1}}\cdots j_0^{a_{i_r}}),
 \quad
 t_M(1;q)\longmapsto\dim M_0,
 \quad
 (\Phi_r,\Psi_r)\longmapsto
 (\Phi_r^{(0)},\Psi_r^{(0)}).
 \label{eq:cusp-replacements}
\end{equation}
\end{theorem}
\begin{proof}
The $\tau \to i\infty$ limit of correlation functions are given by the replacements above. We just need to compare two graph sum rules.

Both sides are summed over partition of $\left[n\right]$ into cycles of at most $1$ type--1 and any number of type--2, i.e. $\cP_n^{\star}$. The relations of cycle coefficients are given \cref{prop:cusp-kernels}
\begin{equation*}
 \frac{\Phi_r^{(0)}(C)}
 {(2\pi\ii)^r\prod_{i\in C}u_i}
 =\Phi_r^{(\mathrm{pl})}(C),
 \quad
 \frac{\Psi_r^{(0)}(C')}
 {(2\pi\ii)^r\prod_{i\in C'}u_i}
 =\Psi_r^{(\mathrm{pl,cyc})}(C')
\end{equation*}
for $r\geq 2$. When $r=1$, our convention gives $\Psi_1^{(0)}=2\pi\ii$ while its plane coefficient is $\Psi_1^{(\mathrm{pl})}(\langle i\rangle)=u_i^{-1}$. Although this gives rise to a term $\Tr_{M_0}j_0^{a_i}$, it has already been shown that this term is zero.
Observe that each $u_i$  occurs exactly once. So the theorem is proved by factoring out $\prod (2\pi \ii u_i)$.
\end{proof}

\appendix
\section{Bell and Bernoulli Polynomials}
\label{app:combinatorics}

This appendix records the two polynomial families used in the coefficient
and constant-term calculations.  We include proofs of the generating
identities, recurrences, and the generalized Bernoulli factorization
needed in \cref{lem:kernel-identities}.

\subsection{Complete exponential Bell polynomials}
\label{app:bell}

\begin{definition}
For $m\geq0$, the complete exponential Bell polynomial is
\begin{equation}
 B_m(x_1,\ldots,x_m)
 :=m!\!\!\sum_{\substack{k_1,\ldots,k_m\geq0\\
 k_1+2k_2+\cdots+mk_m=m}}
 \prod_{j=1}^m\frac1{k_j!}
 \left(\frac{x_j}{j!}\right)^{k_j},
 \qquad B_0:=1.
 \label{eq:Bell-definition}
\end{equation}
\end{definition}

\begin{proposition}[Bell generating function and recurrence]
\label{prop:Bell-properties}
The polynomials in \eqref{eq:Bell-definition} satisfy
\begin{equation}
 \exp\!\left(\sum_{j\geq1}x_j\frac{t^j}{j!}\right)
 =\sum_{m\geq0}B_m(x_1,\ldots,x_m)\frac{t^m}{m!},
 \label{eq:Bell-generating}
\end{equation}
and
\begin{equation}
 B_{m+1}(x_1,\ldots,x_{m+1})
 =\sum_{j=0}^m\binom mj
 B_{m-j}(x_1,\ldots,x_{m-j})x_{j+1}.
 \label{eq:Bell-recurrence}
\end{equation}
\end{proposition}

\begin{proof}
Expanding the exponential as a product gives
\begin{equation*}
 \prod_{j\geq1}
 \left(\sum_{k_j\geq0}\frac1{k_j!}
 \left(\frac{x_jt^j}{j!}\right)^{k_j}\right).
\end{equation*}
The coefficient of $t^m$ is precisely
$B_m/m!$ in \eqref{eq:Bell-definition}, proving
\eqref{eq:Bell-generating}.  Differentiate
\eqref{eq:Bell-generating} with respect to $t$:
\begin{equation*}
 \sum_{m\geq0}B_{m+1}\frac{t^m}{m!}
 =\left(\sum_{j\geq0}x_{j+1}\frac{t^j}{j!}\right)
 \left(\sum_{\ell\geq0}B_\ell\frac{t^\ell}{\ell!}\right).
\end{equation*}
Comparing the coefficient of $t^m$ gives
\eqref{eq:Bell-recurrence}.
\end{proof}

The first values are
\begin{align*}
 B_0&=1,\\
 B_1&=x_1,\\
 B_2&=x_1^2+x_2,\\
 B_3&=x_1^3+3x_1x_2+x_3,\\
 B_4&=x_1^4+6x_1^2x_2+3x_2^2+4x_1x_3+x_4.
\end{align*}
Applying \eqref{eq:Bell-generating} to
$\log H_r=\sum_{j\geq1}L_j\lambda^j/j!$ proves
\eqref{eq:H-Bell} and yields all low-point coefficients in
\cref{sec:low-points}.

\subsection{Generalized Bernoulli polynomials}
\label{app:bernoulli}

\begin{definition}
For an integer $r\geq0$, the generalized Bernoulli polynomials
$B_m^{(r)}(x)$ are defined by
\begin{equation}
 \left(\frac{D}{\e^D-1}\right)^r\e^{xD}
 =\sum_{m\geq0}B_m^{(r)}(x)\frac{D^m}{m!}.
 \label{eq:generalized-Bernoulli}
\end{equation}
For $r=1$ these are the ordinary Bernoulli polynomials, and for $r=0$
one has $B_m^{(0)}(x)=x^m$.
\end{definition}

\begin{proposition}[Elementary properties]
\label{prop:Bernoulli-properties}
For $m,r\geq1$,
\begin{align}
 \partial_xB_m^{(r)}(x)&=mB_{m-1}^{(r)}(x),
 \label{eq:Bernoulli-Appell}\\
 B_m^{(r)}(x+1)-B_m^{(r)}(x)
 &=mB_{m-1}^{(r-1)}(x).
 \label{eq:Bernoulli-difference}
\end{align}
Moreover, $B_m^{(r)}(x)$ is monic of degree $m$.
\end{proposition}

\begin{proof}
Differentiating \eqref{eq:generalized-Bernoulli} with respect to $x$
multiplies its left-hand side by $D$; comparison of coefficients proves
\eqref{eq:Bernoulli-Appell}.  Replacing $x$ by $x+1$ and subtracting
multiplies the generating function by $\e^D-1$, so
\begin{equation*}
 (\e^D-1)\left(\frac{D}{\e^D-1}\right)^r\e^{xD}
 =D\left(\frac{D}{\e^D-1}\right)^{r-1}\e^{xD},
\end{equation*}
which proves \eqref{eq:Bernoulli-difference}.  Finally,
$D/(\e^D-1)=1+O(D)$; hence the coefficient of $x^m$ in
$B_m^{(r)}(x)$ comes solely from $\e^{xD}$ and equals one.
\end{proof}

\begin{proposition}[Factorization at degree $r-1$]
\label{prop:bernoulli-factorization}
For every $r\geq1$,
\begin{equation}
 B_{r-1}^{(r)}(x)=\prod_{j=1}^{r-1}(x-j).
 \label{eq:Bernoulli-factorization}
\end{equation}
\end{proposition}

\begin{proof}
By \cref{prop:Bernoulli-properties}, the left-hand side is monic of
degree $r-1$.  For $j=1,\ldots,r-1$,
\begin{align*}
 \frac{B_{r-1}^{(r)}(j)}{(r-1)!}
 &=[D^{r-1}]
 \left(\frac{D}{\e^D-1}\right)^r\e^{jD}\\
 &=\Res_{D=0}\frac{\e^{jD}}{(\e^D-1)^r}\,\dd D.
\end{align*}
With $t=\e^D$, the last expression is
\begin{equation*}
 \Res_{t=1}\frac{t^{j-1}}{(t-1)^r}\,\dd t=0,
\end{equation*}
because $t^{j-1}$ has degree strictly smaller than $r-1$.  Thus
$1,\ldots,r-1$ are all the roots, and monicity proves
\eqref{eq:Bernoulli-factorization}.
\end{proof}

\end{document}